\documentclass[12pt]{amsart}
\usepackage{amsfonts,amssymb,amscd,amstext,mathrsfs,color, comment}
\usepackage[colorlinks=true, linktocpage,  linkcolor=blue, citecolor=blue, urlcolor=blue]{hyperref}
\input xy
\xyoption{all}

\newtheorem{theorem}{Theorem}[section]
\newtheorem{proposition}[theorem]{Proposition}
\newtheorem{lemma}[theorem]{Lemma}
\newtheorem{corollary}[theorem]{Corollary}    

\theoremstyle{definition}
\newtheorem{definition}[theorem]{Definition}

\newtheorem{remark}[theorem]{Remark}
\numberwithin{equation}{section}

\newcommand{\A}{{\mathscr A}}
\newcommand{\B}{{\mathbb B}}
\newcommand{\C}{{\mathbb C}}
\newcommand{\CC}{{\mathcal C}}
\newcommand{\calA}{\mathcal A}
\newcommand{\calB}{\mathcal B}
\newcommand{\calC}{\mathcal C}
\newcommand{\calW}{\mathcal W}

\newcommand{\NN}{{\mathbb N}}
\renewcommand{\O}{{\mathscr O}}
\newcommand{\F}{{\mathscr F}}
\newcommand{\R}{{\mathbb R}}
\newcommand{\X}{{\mathscr X}}

\newcommand{\Hom}{{\operatorname{Hom}}}

\newcommand{\Ker}{{\operatorname{Ker}\,}}
\newcommand{\inv}{{^{-1}}}
\newcommand{\git}{/\!\!/}
\newcommand{\conj}{{\operatorname{conj}}}

\newcommand{\pt}{\partial}
\renewcommand{\phi}{\varphi}
\newcommand{\GL}{\operatorname{GL}}
\newcommand{\SL}{\operatorname{SL}}
\newcommand{\U}{\operatorname{U}}

\newcommand{\id}{{\operatorname{id}}}
\renewcommand{\int}{\operatorname{int}}
\renewcommand{\Im}{\operatorname{Im}}

\begin{document}

\title{Parametric equivariant Oka principle}

\author{Frank Kutzschebauch, Finnur L\'arusson, Gerald W.~Schwarz}

\address{Frank Kutzschebauch, Institute of Mathematics, University of Bern, Sidlerstrasse 5, CH-3012 Bern, Switzerland}
\email{frank.kutzschebauch@math.unibe.ch}

\address{Finnur L\'arusson, School of Mathematical Sciences, Adelaide University, Adelaide SA 5005, Australia}
\email{finnur.larusson@adelaide.edu.au}

\address{Gerald W.~Schwarz, Department of Mathematics, Brandeis University, Waltham MA 02454-9110, USA}
\email{schwarz@brandeis.edu}

\dedicatory{In memory of Alan T.~Huckleberry}


\date{2 November 2025.  Latest edits 1 September 2026}


\begin{abstract} 

Let $G$ be a reductive complex Lie group and $K$ be a maximal compact subgroup of $G$.  Let $X$ be a reduced Stein $G$-space and $Y$ be a $G$-elliptic manifold.  We prove the following parametric equivariant Oka principle.  The inclusion of the space of holomorphic $G$-maps $X\to Y$ into the space of continuous $K$-maps $X\to Y$ is a weak homotopy equivalence with respect to the compact-open topology.  The proof is divided into a homotopy-theoretic part, which is handled by an abstract theorem of Studer, and an analytic part, for which we prove equivariant versions of the homotopy approximation theorem and the nonlinear splitting lemma that are key tools in Oka theory.  The principle can be strengthened so as to allow interpolation on a $G$-invariant subvariety of $X$ and approximation on a $K$-invariant holomorphically convex compact subset of $X$.
\end{abstract}

\maketitle

\section{Introduction} 
\label{sec:intro}

\noindent
Oka theory is the subfield of complex geometry that is concerned with the homotopy principle in complex analysis.  It has its origin in the pioneering work of Kiyoshi Oka in the late 1930s and was further developed by the Grauert school in the late 1950s through to the early 1970s with a focus on complex Lie groups and homogeneous spaces.  In complex analysis the homotopy principle is known as the Oka principle.  It is an umbrella term for a range of theorems stating that the obstructions to solving various analytic problems on Stein spaces, typically problems that can be cohomologically or homotopically formulated, are purely topological or more precisely homotopy-theoretic in nature.  Oka theory was brought into the modern era in Gromov's seminal paper of 1989 \cite{Gromov1989}, eventually leading to the notions of an Oka manifold, generalising the notion of a homogeneous space, and an Oka map, which are now the central concepts of the theory.  The first major application of Gromov's work was the solution of the Forster conjecture in dimensions greater than 1 \cite{EG1992, Schurmann1997}.  Among the areas in which Oka theory has been applied more recently (with one sample reference for each) are the theory of minimal surfaces \cite{AFLo2021}, the holomorphic Vaserstein problem \cite{IK2012}, complex contact geometry \cite{AFLa2021}, and holomorphic dynamics \cite{AL2022}.  There is an analogous theory in the algebraic category, in some ways similar and in other ways different from analytic Oka theory \cite{LT2019}.  We refer the reader to the monograph \cite{Forstneric2017} and the new survey \cite{Forstneric2025}.

In a series of papers, the authors have brought together Oka theory and geometric invariant theory to develop equivariant Oka theory.  For an overview of this work, see the survey \cite{KLS2022}.  The purpose of the present paper is threefold:
\begin{itemize}
\item  To extend the parametric Oka principle proved in \cite{KLS2018} beyond the setting of homogeneous spaces.
\item  To strengthen the basic Oka principle proved in \cite{KLS2021} to a parametric result.
\item  To combine these two goals in a single theorem proved as simply and cleanly as possible using Studer's abstract framework \cite{Studer2020}.
\end{itemize}
Thus, our main result is the following equivariant parametric Oka principle with interpolation and approximation.

\begin{theorem}   \label{t:main-theorem}
Let $G$ be a reductive complex Lie group and $K$ be a maximal compact subgroup of $G$.  Let $X$ be a reduced Stein $G$-space and $Y$ be a $G$-elliptic manifold.  

{\rm (a)}  The inclusion of the space of holomorphic $G$-maps $X\to Y$ into the space of continuous $K$-maps $X\to Y$ is a weak homotopy equivalence with respect to the compact-open topology.

{\rm (b)}  Let $X'$ be a $G$-invariant subvariety of $X$ and $h:X'\to Y$ be a holomorphic $G$-map.  The inclusion of the space of holomorphic $G$-maps $X\to Y$ that equal $h$ on $X'$ into the space of continuous $K$-maps $X\to Y$ that equal $h$ on $X'$ is a weak homotopy equivalence.

{\rm (c)}  A holomorphic $G$-map $X' \to Y$ extends to a holomorphic $G$-map $X \to Y$ if it extends to a continuous $K$-map $X \to Y$.

{\rm (d)}  Let $L$ be a $K$-invariant holomorphically convex compact subset of $X$, let $X'$ be a $G$-invariant subvariety of $X$, and let $f:X\to Y$ be a continuous $K$-map that is holomorphic on a neighbourhood of $L$ and on $X'$.  Then $f$ can be approximated uniformly on $L$ by holomorphic $G$-maps $X\to Y$ that equal $f$ on $X'$.
\end{theorem}

Part (a) follows from part (b), of course, but is stated separately because until the final section of the paper we focus on (a).  Part (c) is an immediate consequence of (b) and is used to prove (d) in the final section.  The actions of $G$ on $X$ and $Y$ are holomorphic actions by biholomorphisms.  For the definition and basic properties of $G$-ellipticity, see \cite[Section 3]{KLS2021}, where the concept was first defined.  We recall the definition in Section \ref{sec:background}.  Before discussing the proof of the theorem, we list some examples of $G$-elliptic manifolds and cite previous work in which special cases of the theorem were proved.

\begin{remark}   \label{rem:main}
(a)  All $G$-modules and all $G$-homogeneous spaces are $G$-elliptic \cite[Proposition 3.3]{KLS2021}.  In the special case that $Y$ is $G$-homogeneous, Theorem \ref{t:main-theorem}(a) follows from the main theorem of \cite{KLS2018}; see \cite[Theorem E]{KLS2022}.  More generally, the main theorem of \cite{KLS2018} implies Theorem \ref{t:main-theorem}(a) if the $G$-action on $Y$ factors through a transitive action of another complex Lie group, not necessarily reductive, on $Y$.  In Section \ref{sec:Danielewski} we present a class of $G$-elliptic surfaces, most of which are not homogeneous (see (d) below and Remark \ref{rem:not-homogeneous}).

(b)  If $Y$ is a Stein $G$-manifold satisfying the equivariant basic Oka property with jet interpolation ($G$-BOPJI; see \cite{KLS2021}), then $Y$ is easily seen to be $G$-elliptic (see the proof of \cite[Corollary 4.3]{KLS2021}).  Hence, by the main theorem of \cite{KLS2021}, $Y$ is $G$-elliptic if $Y$ is $G$-Oka and all the stabilisers of the $G$-action on $Y$ are finite, in particular if $G$ itself is finite.  (In \cite{KLS2021}, all sources $X$ as in Theorem \ref{t:main-theorem} are taken to be smooth.)  To say that $Y$ is $G$-Oka means that the fixed-point manifold $Y^H$ is Oka for all reductive closed subgroups $H$ of $G$ (see \cite[Section 2]{KLS2021}).  Thus the $G$-Oka property can be investigated using all the resources of non-equivariant Oka theory.

(c)  In view of (a) and (b) it is of interest that there are actions of finite groups on affine spaces that are not known to factor through any transitive action, but with respect to which affine space is equivariantly Oka and hence equivariantly elliptic.  For example, Derksen and Kutzschebauch produced an action of $\C^*$ on $\C^4$ that is not linearisable \cite{DK1998}.  From their construction it is easily seen that the sole nontrivial fixed-point manifold in $\C^4$ is biholomorphic to $\C^2$ and hence Oka.  Thus, $\C^4$ is equivariantly Oka with respect to the $\C^*$-action and therefore equivariantly elliptic with respect to the action of any finite subgroup of $\C^*$.  It seems difficult to determine whether such an action factors through a transitive action of some complex Lie group.

(d)  Danielewski surfaces in $\C^3$ are defined by an equation of the form $xy=f(z)$, where $f$ is an entire function all of whose zeros are simple.  They are $\C^*$-elliptic with respect to the action $t\cdot (x,y,z) = (tx, t^{-1}y, z)$, but most of them are not homogeneous.  Higher-dimensional Danielewski manifolds are hypersurfaces in $\C^n$, $n\geq 4$, defined in a similar way.  Some of them are $\C^*$-elliptic.  A new construction of equivariant sprays with respect to actions of commutative groups (Theorem \ref{thm:G-elliptic}) and other details are given in Section \ref{sec:Danielewski}.

(e)  The famous Koras-Russell cubic 
\[ C=\{ (x,y,s,t)\in\C^4 : x+x^2y+s^2+t^3=0 \} \]
is Oka.  This follows from the result of Leuenberger that $C$ has the volume density property with respect to the volume form $x^{-2}dx\wedge ds\wedge dt$ \cite{Leuenberger2016}.  It is in fact $\C^*$-Oka with respect to the well-known action $\zeta\cdot(x,y,s,t)=(\zeta^6 x, \zeta^{-6}y, \zeta^3 s, \zeta^2 t)$.  Hence, $C$ is $G$-elliptic for every finite subgroup $G$ of $\C^*$.  It is an open question whether $C$ is $\C^*$-elliptic.

To see that $C$ is $\C^*$-Oka, we need to check that the three fixed-point manifolds in $C$, with respect to the subgroups of 2nd, 3rd, and 6th roots of unity, are Oka.  The fixed-point surface $S_2$ of the $\mathbb Z_2$-action is defined by the equations $x+x^2y+t^3=s=0$.  It has the generalised Danielewski surface defined by the equation $u^3v=w^3-1$ in $\C^3$ as a 3-sheeted Galois covering space with covering map $(u, v, w) \mapsto (u^3, v, 0, -uw)$ and covering group generated by the transformation $(u, v, w) \mapsto (\zeta u, v, \zeta^{-1}w)$, where $\zeta$ is a primitive 3rd root of unity.  Similarly, the fixed-point surface $S_3$ of the $\mathbb Z_3$-action is defined by the equations $x+x^2y+s^2=t=0$ and has the generalised Danielewski surface defined by the equation $u^2v=w^2-1$ as a 2-sheeted covering space with covering map $(u, v, w) \mapsto (u^2, v, iuw, 0)$.  It is well known that generalised Danielewski surfaces are Oka -- in fact, they satisfy much stronger algebraic flexibility properties -- so $S_2$ and $S_3$ are Oka.  The fixed-point curve of the $\mathbb Z_6$-action is defined by the equations $x(1+xy)=s=t=0$, so it is the disjoint union of $\C$ and $\C^*$ and hence Oka.  
\end{remark}

To prove Theorem \ref{t:main-theorem}, we make use of the work of Studer \cite{Studer2018, Studer2020, Studer2021}, who developed an abstract framework for proving Oka principles.  His work may be seen as a highly nontrivial adaptation to complex analysis of Gromov's homomorphism theorem \cite[p.~77]{Gromov1986}.  Gromov's theorem states, roughly speaking, that a local weak homotopy equivalence between sheaves of topological spaces is a global weak homotopy equivalence if the sheaves are \emph{flexible}.  Studer's key contribution was to extract from the proofs of some of the fundamental theorems of Oka theory the correct notion of flexibility, allowing him to cleanly separate these proofs into a common abstract homotopy-theoretic part and an analytic part that must be adapted to each particular setting.

Under the hypotheses of Theorem \ref{t:main-theorem}, we let $\pi:X\to Q=X\git G$ be the categorical quotient and define sheaves $\Phi \hookrightarrow \Psi$ on $Q$ by letting $\Phi(U)$, where $U\subset Q$ is open, be the space of holomorphic $G$-maps $\pi^{-1}(U)\to Y$ and $\Psi(U)$ be the space of continuous $K$-maps $\pi^{-1}(U)\to Y$.  With the compact-open topology, these are sheaves of topological spaces, in fact complete metrisable spaces.  (For a summary of the basics on the categorical quotient with references, see the introduction to \cite{KLS2022}.)  By \cite[Theorem 1]{Studer2020}, to conclude that the inclusion $\Phi(Q) \hookrightarrow \Psi(Q)$ is a weak homotopy equivalence and thereby establish Theorem \ref{t:main-theorem}(a), it suffices to prove the following.\begin{itemize}
\item  The inclusion $\Phi\hookrightarrow \Psi$ is a local weak homotopy equivalence.
\item  The quotient $Q$ is covered by open sets $U$ such that every $\mathcal C$-pair $(A,B)$ with $B\subset U$ is weakly flexible for $\Psi$. 
\item  The above property for $\Phi$.
\end{itemize}
We prove the first and second statements, and recall the definitions of weak flexibility and a local weak homotopy equivalence, in Section \ref{sec:local}.  The proofs do not require the ellipticity assumption on $Y$.  The bulk of the paper is devoted to the proof of the third statement.  The proof is presented in Section \ref{sec:hol-weak-flex}, using the equivariant parametric homotopy approximation theorem proved in Section \ref{sec:approximation} (Theorem \ref{thm:K-ell-implies-EPHAP}) and the equivariant nonlinear splitting lemma proved in Section \ref{sec:splitting} (Theorem \ref{thm:main-splitting}).  The former is an equivariant version of \cite[Theorem 6.6.2]{Forstneric2017}, a key tool in Oka theory, and the latter is an equivariant version of a consequence of another key tool, \cite[Proposition 5.8.4]{Forstneric2017}.  In the final section we show how interpolation can be incorporated into the proof of Theorem \ref{t:main-theorem}(a) so as to prove Theorem \ref{t:main-theorem}(b) and (d).

 \smallskip
 
\noindent  {\bf Acknowledgements.}  F.~Kutzschebauch was supported by Schweizerischer Nationalfonds grants 200021-207335 and 2000-1-240004.  F.~L\'arusson and G.~W.~Schwarz thank the University of Bern, where much of this work was done, for its hospitality.  The authors thank a referee for pointing out a gap in the proof of a previous version of Theorem \ref{thm:main-splitting}.  They also thank Song-Yan Xie for help with Remark \ref{rem:main}(e).

\smallskip

\noindent  {\bf Use of artificial intelligence.}  In the final stage of preparation of this paper, the AI assistant Claude was used to weed out minor mistakes, inconsistent notation, and typographical errors.

\section{Background and preparation} 
\label{sec:background}

\subsection{Equivariant ellipticity}  
A manifold $Y$ is said to be elliptic if it carries a dominating spray, that is, there is a holomorphic map $s:E\to Y$, called a spray, defined on the total space of a holomorphic vector bundle $E$ on $Y$, such that $s(0_y)=y$ for all $y\in Y$, which is dominating in the sense that $s\vert_{E_y}:E_y\to Y$ is a submersion at $0_y$ for all $y\in Y$.  Suppose that a complex Lie group $G$ acts on $Y$.  (Such an action is always assumed to be holomorphic.)  We say that $s$ is a $G$-spray if the action on $Y$ lifts to an action on $E$ by vector bundle isomorphisms such that both $s$ and the projection $E\to Y$ are equivariant.  We say that $Y$ is $G$-elliptic if it carries a dominating $G$-spray.  This notion was introduced in \cite[Section 3]{KLS2021}.  Similarly, we define $K$-ellipticity of $Y$ for a real Lie group $K$ acting continuously and hence real-analytically on $Y$ by biholomorphisms.

\begin{proposition}   \label{prop:gell.eq.kell} 
Let $G$ be a reductive complex Lie group, $K$ be a maximal compact subgroup of $G$, and $Y$ be a $G$-manifold.  Then $Y$ is $G$-elliptic if and only if it is $K$-elliptic.
\end{proposition}

\begin{proof}
Clearly, if $Y$ is $G$-elliptic, then it is $K$-elliptic.  Conversely, suppose that $Y$ is $K$-elliptic and that $\sigma: E\to Y$ is a $K$-equivariant dominating spray, where $E$ is a holomorphic $K$-vector bundle on $Y$.  By \cite[\S 6, Proposition 1]{HK1995}, $E$ is naturally a $G$-vector bundle and since $\sigma: E\to Y$ is holomorphic and $K$-equivariant, it is $G$-equivariant.  Hence $Y$ is $G$-elliptic.
\end{proof}

\subsection{Equivariant embeddings}
Let $G$ be a reductive complex Lie group, $K$ be a maximal compact subgroup of $G$, and $X$ be a Stein $G$-space, here and throughout assumed to be reduced.  Let $\pi:X\to X\git G$ be the categorical quotient map.  The following proposition will be used later on.

\begin{proposition}  \label{prop:equivariant-embedding}
{\rm (a)}  If $U$ is a relatively compact Stein open subset of $X\git G$, then the Stein open subset $\pi^{-1}(U)$ of $X$ has a $G$-equivariant proper holomorphic embedding into a $G$-module.

{\rm (b)}  If $V$ is a relatively compact $K$-invariant Stein open subset of $X$, then $V$ embeds $K$-equivariantly as an open subset of a $G$-subvariety in a $G$-module.
\end{proposition}

\begin{proof}
{\rm (a)}  By Heinzner's embedding theorem \cite{Heinzner1988}, it suffices to show that $\pi^{-1}(U)$ contains only finitely many slice types.  The slice type of a point $x\in X$, whose orbit $Gx$ is closed, is the isomorphism class of the isotropy representation of $G_x$ on $T_x X$.  By the holomorphic slice theorem, a sufficiently small saturated neighbourhood of any closed orbit $Gx$ in $X$ contains only finitely many slice types.  Since $U$ is relatively compact, $\pi^{-1}(U)$ contains only finitely many slice types.

{\rm (b)}  Apply (a) to a relatively compact Stein neighbourhood $U$ of $\pi(V)$.
\end{proof}

\subsection{Stein compact sets and Kempf-Ness sets}   \label{sec:Stein.compacts}
Let $G$, $K$, and $X$ be as above.  For the following, see  \cite[p.~341]{HK1995}.  There is a real-analytic $K$-invariant strictly plurisubharmonic exhaustion function $\phi: X\to [0,\infty)$ and an associated real-analytic subvariety $R$ of $X$, called a Kempf-Ness set, with the following properties.
\begin{itemize}
\item  $R$ consists of precisely one $K$-orbit in every closed $G$-orbit in $X$.
\item  The inclusion $R\hookrightarrow X$ induces a homeomorphism $R/K\to X\git G$, where the orbit space $R/K$ carries the quotient topology. 
\item  For every neighbourhood $U$ of $R$, we have $G\cdot U=X$.
\item  When $X$ is smooth, $R$ is a $K$-equivariant continuous strong deformation retract of $X$, such that the deformation preserves the closure of each $G$-orbit \cite{HH1994}.
\end{itemize}
For $c>0$, let $X_c:=\phi\inv([0,c))$.  Note that $X_c$ is $K$-stable\footnote{We use the synonyms {\it stable} and {\it invariant} interchangeably.} and is the interior of $\overline{X_c}=\phi\inv([0,c])$.

\begin{proposition}   \label{prop:HK} 
\begin{enumerate}
\item For any $c>0$, $\overline{X_c}$ is $\O(X)$-convex.
\item For any $c>0$, $X_c$ is Stein and Runge in $X$.
\end{enumerate}
\end{proposition}

\begin{proof}
By \cite[Theorem 2.5.2]{Forstneric2017}, we have (1).  For (2), if $M\subset X_c$ is compact, then it is contained in some $\overline{X_{c'}}$ for $0<c'<c$.  The $\O(X_c)$-convex hull of $M$ is contained in the $\O(X)$-convex hull of $M$ which is a compact subset of $\overline{X_{c'}}\subset X_c$.  Thus $X_c$ is holomorphically convex and open in $X$, hence Stein.  If $f\in\O(X_c)$, then its restriction to any $\overline{X_{c'}}$, $0<c'<c$, is uniformly approximable by elements of $\O(X)$.  Hence $X_c$ is Runge in $X$.
\end{proof}

\begin{lemma}   \label{lem:intersection.Stein}
Let $\Omega$ be a Stein open set in the complex $K$-space $Z$. Then 
\[  \Omega':=\bigcap_{k\in K} k\cdot \Omega \]
is open, $K$-invariant, and Stein.
\end{lemma}

\begin{proof}
For each $k\in K$, $k\cdot \Omega$ is Stein.  Since $K$ is compact, $K\cdot(Z\setminus\Omega)$ is closed in $Z$, hence its complement $\Omega'$ is open.  Thus $\Omega'$ is Stein if it is holomorphically convex.  Let $M\subset\Omega'$ be compact.  The $\O(\Omega')$-convex hull $\widehat M$ of $M$ is contained in the (compact) $\O(k\cdot\Omega)$-convex hull of $M$ for all $k$.  Hence, $\widehat M$ is a compact subset of $\Omega'$ and $\Omega'$ is Stein.  
\end{proof}

Using \cite{Siu1976},  we obtain the following.

\begin{corollary}   \label{cor:Stein-nbhd}
Let $M$ be a closed Stein $K$-stable subspace of the complex $K$-space $Z$.  Then any neighbourhood of $M$ in $Z$ contains a neighbourhood which is $K$-invariant and Stein.
\end{corollary}

\subsection{Sprays and parametric sprays}  \label{sec:sprays}  
The results in this section are used in Sections \ref{sec:approximation}, \ref{sec:splitting}, and \ref{sec:hol-weak-flex}.  Let $Y$  be $G$-elliptic with corresponding $G$-vector bundle $E$ and dominating $G$-equivariant spray   $s: E\to Y$.  Let $E''_y=\Ker(Ds)_0: E_y\to T_yY$ for $y\in Y$.  Since $s$ is dominating, $E''$ is a $G$-vector subbundle of $E$ and $(Ds)_0$ induces a $G$-isomorphism of $E':=E/E''$ and $TY$.

Let $X$ be a Stein $G$-space as before and let $f: X\to Y$ be a $G$-equivariant holomorphic map.  Let $F=f^*E$ and $\sigma=f^*s: F\to Y$.  Then $\sigma\vert_{F_x}=s\vert_{E_{f(x)}}$, so $\sigma$ is dominating and $G$-equivariant with core $f$ (meaning that $\sigma=f$ on the zero section of $F$).  Since $X$ is Stein, we have the following result (see \cite[Lemma 7.2]{KLS2021} for some basic facts about equivariant vector bundles on a Stein space).
 
\begin{lemma}   \label{lem:structure.E'}
Let $Y$, $f$, etc.\ be as above.  Let $F''=f^*E''$.  Then $F''$ admits a complementary $G$-vector subbundle $F'$ of $F$ and $D(\sigma\vert_{F'})_0: F'\to f^*TY$ is a $G$-isomorphism. 
\end{lemma}

Let $f$, $F$ and $\sigma$ be as above.  Let $\gamma_f$ denote the function $x\mapsto (x,f(x))$, $x\in X$.  Let $\rho: F\to X$ be the bundle projection and let $\Gamma(F)$ denote the holomorphic sections of $F$.  If $\xi\in\Gamma(F)$, let $\Im\xi$ denote  its image in $F$ and let $\Xi$ denote the image of the zero section.

\begin{lemma}   \label{lem:local.biholom} 
Let $f$, $F$, etc.\ be as above. Assume that $F''$ is the zero bundle.
\begin{enumerate}
\item There is a Stein $K$-neighbourhood $\Omega$ of $\Xi$ such that the map
\[  \Phi: \Omega\to X\times Y,\quad v\mapsto (\rho(v),\sigma(v))  \]
is $K$-equivariant and biholomorphic onto its (open) image.
\item  If $\xi\in\Gamma(F)$ with $\Im\xi\subset\Omega$, then $\Phi(\Im\xi) = \gamma_{f'}(X)$ where $f': X\to Y$ is holomorphic.  Conversely, if $f': X\to Y$ is   holomorphic and $\gamma_{f'}(X)\subset\Phi(\Omega)$, then $\Phi\inv(\gamma_{f'}(X)) = \Im\xi$, where $\xi\in\Gamma(F)$.
Moreover, $\xi$ is $K$-equivariant if and only if $f'$ is $K$-equivariant. 
\end{enumerate}
\end{lemma}

\begin{proof}
Choose a $K$-invariant norm $\vert\cdot\vert$ on $F$.  For any $x\in X$, there is $\epsilon>0$ such that $\Phi\vert_{F\vert_{Kx}}$ is a $K$-biholomorphism from $\{\xi\in F\vert_{Kx}: \vert \xi\vert<\epsilon\}$ onto a $K$-neighbourhood of $\gamma_f(x)$ in $X\times Y$.  Clearly for $x'$ sufficiently close to $x$, $\Phi\vert_{F\vert_{Kx'}}$ is a $K$-biholomorphism from $\{\xi\in F\vert_{Kx'}: \vert \xi\vert<\epsilon/2\}$ onto a $K$-neighbourhood of $\gamma_f(x')$ in $X\times Y$.  Thus there is a neighbourhood $\Omega$ of the Stein subset $\Xi\subset F$ on which $\Phi$ is a $K$-biholomorphism.  By Corollary \ref{cor:Stein-nbhd}, we may assume that $\Omega$ is $K$-stable and Stein.
\end{proof}

\begin{remark}   \label{rem:mult.by.t}
Using a $K$-invariant strictly plurisubharmonic function $\phi$ as in \cite[Proposition 3.3.1]{Forstneric2017}, we may arrange that the fibres $\Omega_x$ are convex.
\end{remark}

We now consider parametric sprays.  Let $P$ be a compact Hausdorff space and let $f: X\times P\to Y$ be continuous, $G$-equivariant, and holomorphic for each fixed $p\in P$.  We assume that $P$ is a finite polyhedron, so $P\subset\R^n\subset\C^n$ for some $n$.  Let $Z=\C^n\times X\times Y$.  Let $L$ be a $K$-stable $\O(X)$-convex compact subset of $X$. Let $\gamma_f: P\times X\to Z$ send $(p,x)$ to $(p,x,f(p,x))$ and set $M:=\gamma_f(P\times L)$. 

\begin{lemma}   \label{lem:U.exists}
There is a $K$-invariant Stein neighbourhood $U$ of $M$ in $Z$.
\end{lemma}

\begin{proof}
By \cite[Corollary 3.6.6]{Forstneric2017}, there is a Stein neighbourhood $U$ of $M$ in $Z$ which by Lemma \ref{lem:intersection.Stein} we may assume is $K$-invariant.  
\end{proof}

Let $\pi_Y: U\to Y$ be the projection.  Then $F:=\pi_Y^*E$ is a holomorphic $K$-vector bundle over the Stein $K$-space $U$.  Moreover, $\sigma =\pi_Y^*s$ is a dominating spray map with core $\pi_Y$.  Since $P$ is compact, there is a neighbourhood $V\Subset X$ of $L$ such that $\gamma_f(P\times V)\Subset U$.  Since $L$ is $\O(X)$-convex and $K$-stable, we may assume that $V$ is Stein and $K$-stable.  Let $F''$ denote the kernel of $(D\sigma)_0\subset F$.

\begin{lemma}   \label{lem:parametric-splitting}
Let $f$, $L$, $V$, etc.\ be as above.  Let  $\tilde F$ and $\tilde F''$ denote the restrictions of  $F$ and $F''$ to $P\times V$.
\begin{enumerate}
\item There is a continuous family $\tilde F'_p$ of holomorphic $K$-subbundles  of $\tilde F_p$ which are complementary to $\tilde F''_p$, $p\in P$.
\item The splittings of $0\to \tilde F'' \to \tilde F $ correspond to continuous families of holomorphic $K$-equivariant sections of $\Hom(\tilde F'_p,\tilde F''_p)$, $p\in P$.
\item $\tilde F'$ and $\tilde f^*TY$ are isomorphic as holomorphic $K$-vector bundles.
\end{enumerate}
\end{lemma}

\begin{proof}
Since $U$ is Stein, there is a $K$-subbundle $F'$ of $F$ complementary to $F''$. Now use \cite[Lemma 7.2]{KLS2021} and restrict to $P\times V$.
\end{proof}

The following may not be necessary, but it is enough to get what we eventually need.  We  add the assumption that $P$ is contractible, so there is a deformation retraction of $P$ to a point $p_0\in P$.

\begin{lemma}  \label{lem:deform}
Over $P$ we have  a continuous family of $K$-equivariant holomorphic bundle isomorphisms $P\times \tilde F_{p_0}\simeq \tilde F$ and $P\times\tilde F'_{p_0}\simeq\tilde F'$.
\end{lemma}

\begin{proof}
Let $h: I\times P\to P$ be the deformation retraction and let $\tilde h$ be the map
\[  I\times P\times V\to I\times Y, \quad (t,p,x)\mapsto (t,f_{h(t,p)}(x)).  \]
Let $\tilde E=\tilde h^*(I\times E)$.  As in \cite[Theorem 3.8]{KLS2018}, we have $\tilde E_0\simeq \tilde E_1$.  But $\tilde E_1\simeq P\times \tilde F_{p_0}$ while $\tilde E_0\simeq \tilde F$. The same argument works for $\tilde F'$.
\end{proof}

Note that the fibre dimension of $F'$ is $\dim Y$.  Let $\pi: F'\to V\times P$ be the bundle projection and let $Z=V\times\C^n\times Y$ where $P\subset\R^n\subset\C^n$. Let $\Gamma(F')$ denote the holomorphic $P$-families of $F'$, they are continuous sections which are holomorphic on each $\{p\}\times V$.

\begin{theorem}\label{thm:eq-biholom}  
Let $\Theta\simeq V\times P$ denote the zero section of $F'$.
\begin{enumerate}
\item  There is a Stein $K$-neighbourhood $\Omega$ of $\Theta$ in $Z$ such that 
\[ \Phi:  \Omega\to Z,\quad v\mapsto (\pi(v),\sigma(v)), \]
is continuous and $K$-equivariant such that each $\Phi_p: \Omega_p\to V \times Y$ is $K$-biholomorphic onto its (open) image.
\item  If $\psi\in\Gamma(F')$   with $\Im\psi\subset \Omega$, then
$\Phi(\Im\psi)=\gamma_{f'}(V\times P)$ where  $f': V\times P\to Y$ is a  holomorphic $P$-family. Conversely, if $f': V\times P  \to Y$ is a  holomorphic $P$-family and $\gamma_{f'}(V\times P)\subset\Phi(\Omega)$, then $\Phi\inv(\gamma_{f'}(V\times P))=\Im\psi$ where $\psi\in\Gamma(F')$.
Moreover, $\psi$ is $K$-equivariant if and only if $f'$ is $K$-equivariant.
\end{enumerate}
\end{theorem}

\begin{proof}
Part (2) follows from part (1), which is proved exactly as in Lemma \ref{lem:local.biholom}.
\end{proof}

\section{Danielewski manifolds}
\label{sec:Danielewski}

\subsection{Sufficient condition for $G$-ellipticity}
Let $X$ be a complex manifold with the action of a reductive complex Lie group $G$.  Let $\A(X)$ denote the holomorphic vector fields on $X$ and let $\mathscr X(G)$ denote the character group of $G$.  We have an action of $G$ on $\A(X)$,
$$
G\times\A(X)\times X\ni (g,\xi,x)\mapsto (g_*\xi)(x)=Dg\vert_{g\inv x}(\xi(g\inv x)).
$$
Alternatively,  $(g_*\xi)(f)=(\xi(f\circ g))\circ g\inv$ for $f\in\O(X)$.

\begin{remark}
A calculation shows that for $g$, $h\in G$, $(gh)_*=g_*\circ h_*$.
\end{remark}

Let $\chi\in\mathscr X(G)$. We say that $G$ acts on $\xi\in\A(X)$ by $\chi$ and write that $\xi\in\A(X)_\chi$ if $g_*\xi=\chi(g)\xi$, $g\in G$.  Let $\conj(g)$ denote the conjugation action of $g$ on $G$. 

\begin{remark}   \label{rem:G.action}
Suppose that $G$ is commutative. Then for any character $\chi$ of $G^0$ and $g\in G$, $\chi_g=\chi\circ\conj(g\inv)=\chi$ so that $G$ preserves $\A(X)_\chi$.
\end{remark}

\begin{theorem}   \label{thm:G-elliptic}
Let $X$ be a complex $G$-manifold where $G$ is commutative. Assume that there are   $\chi_1,\dots,\chi_n\in\X(G)$ such that finitely many complete elements of the $\A(X)_{\chi_j}$ span the tangent space of $X$ at every point. Then $X$ is $G$-elliptic.
\end{theorem}

\begin{proof}
By hypothesis there are complete vector fields $\xi_{(i,1)},\dots,\xi_{(i,m_i)}\in\A(X)_{\chi_i}$, $i=1,\dots,n$, that span the tangent space of $X$ at every point.  Let $k=\sum_im_i$.  Let $\phi_{(i,j)}^s$ denote the flow of $\xi_{(i,j)}$, $j=1,\dots,m_i$, $i=1,\dots,n$.  Define $\phi_{(a_1,\dots,a_k)}: X\to X$ by
$$
x\mapsto(\phi_{(1,1)}^{a_1}\circ\dots\circ\phi_{(1,m_1)}^{a_{m_1}}\circ\dots\circ\phi_{(n,m_n)}^{a_k})(x).
$$
We view $\phi$ as a spray map on the trivial bundle $\C^k\times X$ with image in $X$. Then  
$$
(g\circ\phi_{(a_1,\dots,a_k)}\circ g\inv)(x)=\phi_{(\chi_1(g)a_1,\dots,\chi_n(g)a_k)}(x), \quad g\in G,\ x\in X.
$$
Now let $G$ act on the basis vector $e_{(i,j)}\in\C^k$ by $\chi_i$. Then with this new action on $\C^k$, which we now call $V$, we get a  dominating spray   $\psi: V\times X\to X$ which is $G$-equivariant.
\end{proof}

\begin{remark}
The proof above produces local equivariant sprays even when the vector fields are not complete.  This does not work for a non-commutative group: the spray given by composition of local flows of equivariant vector fields need not be equivariant.  Local sprays produced from local flows of vector fields are a key tool in standard Oka theory, but are usually not available in the equivariant case.  This is the reason we require $G$-ellipticity in the proof of Theorem \ref{t:local-weak-eq-with-interpolation}.
\end{remark}

\subsection{Danielewski manifolds}
Let $p:\C^n\to\C$ be a holomorphic function whose zero set is smooth and reduced. That is, if $p( x)=0$, then at least one of the partial derivatives $\pt p/\pt x_i(x)$ does not vanish. Let 
$$
X=D_p:=\{(u,v, x) : uv-p( x)=0\}\subset \C^{n+2}.
$$
It is easily seen that $X$ is smooth of dimension $n+1$.  As shown in \cite{KK2008}, $X$ has the density property and is therefore elliptic.  We have an action of $T=\C^*$ on $\C^2$ by $t\cdot (a,b) = (ta,t\inv b)$, which extends by the trivial action on $\C^n$ to an action on $X$.  Let $(u,v)$ be the corresponding coordinate functions.  Note that the natural action of $T$ on functions on $X$ is via $f\mapsto f\circ t\inv$.  Then $t\cdot u=u\circ t\inv = t\inv u$ and $t\cdot v=tv$.

If $X$ is $T$-elliptic, then $X^T\simeq \{ x\in\C^n : p(x)=0\}$ is elliptic \cite[proof of Proposition 3.2]{KLS2021}.  To obtain a converse we need to assume more.

\begin{proposition}  \label{prop:main}
Suppose that   there are complete vector fields $\xi_1,\dots,\xi_m$ on $\C^n$ with the following property.  The $\xi_j$ annihilate $p$ and their restrictions to $X^T$ generate $\A(X^T)$ as an $\O(X^T)$-module. Then $X$ is $T$-elliptic.
\end{proposition}

We will apply Theorem \ref{thm:G-elliptic}.  First we need some preliminaries.  For $i=1,\dots,n$, let 
\begin{align*}
\nu_i&=u\frac {\pt}{\pt x_i}+\frac{\pt p}{\pt x_i}\frac{\pt}{\pt v},\\
\nu_i'&=v\frac {\pt}{\pt x_i}+\frac{\pt p}{\pt x_i}\frac{\pt}{\pt u},\\
H&=u\frac{\pt}{\pt u}-v\frac{\pt}{\pt v}.
\end{align*}
These vector fields annihilate $uv-p( x)$, hence can be considered as vector fields on $X$.  Let $$\Delta_{ij}=\frac{\pt p}{\pt x_i}\frac{\pt}{\pt x_j}-\frac{\pt p}{\pt x_j}\frac{\pt}{\pt x_i}.$$
We leave the proofs of the following lemmas to the reader.

\begin{lemma}\label{lem:trivial}
\begin{enumerate}
\item $[\nu_i,\nu_j]=[\nu_i',\nu_j']=0$ for all $i$, $j$.
\item $\nu_1,\ldots,\nu_n$ are complete holomorphic vector fields of weight $-1$ and $\nu_1',\ldots,\nu_n'$ are complete holomorphic vector fields of weight $1$.  If $p$ is a polynomial, then $\nu_i$ and $\nu_i'$ are all LNDs.
\item $H$ is complete of weight $0$.
\item For $i<j$,
$$[\nu_i,\nu_j']=\frac{\pt^2 p}{\pt x_i\pt x_j}H+\Delta_{ij},$$
which is a vector field of weight $0$.
\end{enumerate}
\end{lemma}

When $u\neq 0$, the projections of the $\nu_i$ to $\C^n$ are linearly independent and $H\neq 0$.  Hence, the $\nu_i$ and $H$ span $TX$.  A similar result holds if $v\neq 0$.  Thus we only need to worry about the case that $u=v=0$, that is, when $x\in X^T$.

\begin{lemma}
Suppose that $x\in X^T$.  
\begin{enumerate}
\item The span of  the $\nu_i$ and $\nu_i'$ at $ x$ is that of $\dfrac{\pt}{\pt u}$ and $\dfrac{\pt}{\pt v}$.
\item The span of the $\Delta_{ij}$ is an $(n-1)$-dimensional subspace of $\C^n$, $i\neq j$.
\end{enumerate}
\end{lemma}

\begin{proof}[Proof of Proposition \ref{prop:main}]
The hypotheses of Theorem \ref{thm:G-elliptic} would be satisfied if the $[\nu_i,\nu_j']$ were complete vector fields, but this we cannot assert.  We are saved by the vector fields $\xi_j\in\A(\C^n)$.  They extend to complete vector fields on $\C^{n+2}$ which annihilate $uv-p( x)$ and by hypothesis their restrictions to $X^T$ generate $\A(X^T)$ over $\O(X^T)$.
\end{proof}

\begin{corollary}\label{cor:n=1}
If $n=1$, then $X^T$ consists of isolated reduced points and hence $X$ is $T$-elliptic.
\end{corollary}

\begin{remark}  \label{rem:not-homogeneous}
Only a few Danielewski surfaces are homogeneous with respect to an action of a complex Lie group.  If $p \in \O(\C)$ has exactly one zero, then $D_p$ is $T$-biholomorphic to the 2-dimensional representation with weights 1 and $-1$. If $p \in \O(\C)$ has exactly two zeros, then $D_p$ is $T$-biholomorphic to the affine quadric $\SL_2 (\C)/H$, where $H$ is the maximal torus and $T$ acts by left multiplication.  The $T$-ellipticity in those cases was established in our earlier paper \cite{KLS2021}.  If, however, $p$ has more than two zeros (possibly infinitely many), then $D_p$ cannot be a homogeneous space of a complex Lie group. Indeed, $D_p$ has trivial fundamental group and strongly deformation-retracts onto a chain of at least two spheres \cite[Section 3.3]{Lind2006}, so it is not on the list of complex homogeneous surfaces in \cite{Huckleberry1986}.
\end{remark}

\begin{proposition}  \label{prop:n=2}
If $n\ge 2$ and $p\in \C[x_1, x_2, \ldots ,x_n]$ is a polynomial which is linear in each variable separately and whose zero set is smooth and reduced, then $X$ is $T$-elliptic.
\end{proposition}
 
\begin{proof}
The vector fields $\Delta_{ij}$ on $\C^n$ annihilating $p$ are complete and span the tangent space of $X^T \simeq \{ x\in\C^n : p(x)=0\}$ at every point \cite[Lemmas 5.2 and 5.3]{IK2012}.  Thus Proposition \ref{prop:main} applies.
\end{proof}

The same proof shows that $X= \{(u,v, x)\in\C^{n+2} : uv-g(p(x)) =0\}$ is $T$-elliptic for any $p\in \C[x_1, \ldots ,x_n]$ that is linear in each variable separately and any $g\in \O(\C)$ such that $g \circ p$ has a smooth and reduced zero set.  If $g$ has infinitely many zeros (but is not identically zero), then the second homology of $X$ is not finitely generated \cite[Proposition 4.1]{KZ1999}, so $X$ is not biholomorphic to a homogeneous space of a complex Lie group. 

\begin{remark}  
We have seen that $X^T$ being elliptic is necessary for $X$ to be $T$-elliptic.  Of course, when $n\geq 2$, $X^T$ need not be elliptic.  Conversely, if $X^T$ is elliptic (or, equivalently, Oka, as $X^T$ is Stein), then $X$ is $T$-Oka, therefore $H$-Oka for every finite subgroup $H$ of $T$, and hence $H$-elliptic (see Remark \ref{rem:main}(b)).  It is an interesting open question whether $H$-ellipticity for every finite subgroup $H$ of $T$ implies $T$-ellipticity.
\end{remark}

\section{Topological flexibility and local weak homotopy equivalence} 
\label{sec:local}

\noindent
We begin by recalling key definitions from \cite[Section 1.2]{Studer2020}.  We denote the closed unit ball in $\R^n$, $n\geq 0$, by $\B_n$ and its boundary by $\partial\B_n$.  We take $\B_n$ to be a point and $\partial \B_n$ to be empty when $n=0$.  Also, write $I=[0,1]$.

Let $\Phi$ and $\Psi$ be sheaves of topological spaces over a topological space $Q$.  A morphism $\alpha:\Phi\to\Psi$ is said to be a \emph{local weak homotopy equivalence} if whenever $U$ is a neighbourhood of a point $p$ in $Q$ and $f:\B_n\to\Psi(U)$ is a continuous map whose restriction to $\partial\B_n$ factors through $\alpha_U$ by a continuous map $\phi:\partial \B_n\to\Phi(U)$, there is a neighbourhood $V\subset U$ of $p$ such that in the commuting square below, $\rho\circ f$ can be deformed, keeping the square commuting, until there is a lifting in the square.  Here, both restriction maps $\Phi(U)\to \Phi(V)$ and $\Psi(U)\to\Psi(V)$ are denoted by $\rho$.
\[ \xymatrix{
\partial \B_n \ar^{\rho\circ\phi}[r] \ar_j[d] & \Phi(V) \ar^{\alpha_V}[d]  \\ \B_n \ar_{\rho\circ f}[r] \ar@{-->}[ur] & \Psi(V)
} \]

It is convenient to have the following lemma.

\begin{lemma}   \label{l:homotopy-equiv}
Suppose that every point in $Q$ has arbitrarily small neighbourhoods $U$ such that the induced map $\alpha_U:\Phi(U)\to\Psi(U)$ is a weak homotopy equivalence.  Then $\alpha$ is a local weak homotopy equivalence.
\end{lemma}

\begin{proof}
Let $U$ and $f$ be as above.  We may assume that $\alpha_U$ is a weak homotopy equivalence.  We will verify the defining property above with $V=U$.  Since the inclusion $j:\partial \B_n\hookrightarrow \B_n$ is a cofibration, the precomposition maps
\[ j_\Phi^*:\mathscr C(\B_n,\Phi(U))\to\mathscr C(\partial \B_n,\Phi(U)), \qquad j_\Psi^*:\mathscr C(\B_n,\Psi(U))\to\mathscr C(\partial \B_n,\Psi(U)) \]
are Hurewicz fibrations (hence Serre fibrations).  Since $\alpha_U$ is a weak homotopy equivalence, the postcomposition maps
\[ {\alpha_U}_*:\mathscr C(\B_n,\Phi(U))\to\mathscr C(\B_n,\Psi(U)),  \qquad  {\alpha_U}_*:\mathscr C(\partial \B_n,\Phi(U))\to\mathscr C(\partial \B_n, \Psi(U))  \]
are weak homotopy equivalences.  Consider the fibres $F_\Phi=(j_\Phi^*)^{-1}(\phi)$ and $F_\Psi=(j_\Psi^*)^{-1}(f\circ j)$.  By the long exact sequence of homotopy groups for a Serre fibration, the map ${\alpha_U}_*:F_\Phi\to F_\Psi$ is a weak homotopy equivalence; in particular it induces a surjection of path components.  Hence, $f\in F_\Psi$ can be deformed within $F_\Psi$ to a map in ${\alpha_U}_*(F_\Phi)$, as desired.
\end{proof}

Next we recall the definition of \emph{weak flexibility} for $\Psi$ of a pair $(A, B)$ of compact subsets of $Q$.  Let $U$, $V$, and $W$ be neighbourhoods of $A$, $B$, and $A\cap B$, respectively, and $a:\B_n\to\Psi(U)$, $b:\B_n\to\Psi(V)$, and $c:\B_n\times I\to\Psi(W)$ be continuous maps such that $a\vert_W=c_0$, $b\vert_W=c_1$, and $c_s\vert_{\partial\B_n} = c(\cdot, s)\vert_{\partial\B_n}$ is independent of $s\in I$.  Then there are smaller neighbourhoods $U'$ of $A$, $V'$ of $B$, and $W'$ of $A\cap B$, and homotopies $a_t:\B_n\to\Psi(U')$, $b_t:\B_n\to\Psi(V')$, and $c_{s,t}:\B_n\to\Psi(W')$ with $a_0=a\vert_{U'}$, $b_0=b\vert_{V'}$, and $c_{s,0}=c_s\vert_{W'}$, such that:
\begin{itemize}
\item  $c_{0,t}=a_t\vert_{W'}$ and $c_{1,t}=b_t\vert_{W'}$ for all $t\in I$,
\item  $a_t\vert_{\partial\B_n}$, $b_t\vert_{\partial\B_n}$, and $c_{s,t}\vert_{\partial\B_n}$ are independent of $t\in I$,
\item  $c_{s,1}$ is independent of $s\in I$, so $a_1\vert_{W'}=b_1\vert_{W'}$,
\item  $a_t\vert_{A^\circ}$ is in a prescribed neighbourhood of $a_0\vert_{A^\circ}:\B_n\to\Psi(A^\circ)$ with respect to the compact open topology, for all $t\in I$.  Here, $A^\circ$ denotes the interior of $A$.
\end{itemize}

We now turn to the proof of Theorem \ref{t:main-theorem}(a).  As before, we let $G$ be a reductive complex Lie group, $K$ be a maximal compact subgroup of $G$, $X$ be a Stein $G$-space, $\pi:X\to Q=X\git G$ be the categorical quotient, and $Y$ be a $G$-manifold.  The results in this section do not require $Y$ to be $G$-elliptic.  We recall that the sheaves $\Phi \hookrightarrow \Psi$ on $Q$ are defined by letting $\Phi(U)$, where $U\subset Q$ is open, be the space of holomorphic $G$-maps $\pi^{-1}(U)\to Y$ and $\Psi(U)$ be the space of continuous $K$-maps $\pi^{-1}(U)\to Y$ with the compact-open topology. 

We begin with the easiest of the three parts of the proof of Theorem \ref{t:main-theorem}(a).

\begin{proposition}   \label{p:top-flex}
The quotient $Q$ is covered by open sets $U$ such that every $\mathcal C$-pair $(A,B)$ with $B\subset U$ is weakly flexible for $\Psi$. 
\end{proposition}

The notion of a $\mathcal C$-pair is defined below (Definition \ref{C-pairs}), but the proof only requires $A$ and $B$ to be compact subsets of $Q$.

\begin{proof}
We verify the stronger flexibility property introduced and applied by Gromov in \cite[Sections 1.4.2 and 2.2.1]{Gromov1986}.  It does not allow the map $a:\B_n\to\Psi(U)$ above to be deformed, that is, the homotopy $a_t$ is required to be constant.

Take any compact subsets $A$ and $B$ of $Q$ and let $C=A\cap B$.  Let $V$ and $W$ be neighbourhoods of $B$ and $C$, respectively, with $W\subset V$.  Let $b:\B_n\to\Psi(V)$ be continuous and $c:\B_n\times I\to\Psi(W)$ be a homotopy with $c(\cdot, 1) = b\vert_W$.  Then the restriction of $c$ to a smaller neighbourhood of $C$ extends to a homotopy $\tilde b:\B_n\times I\to\Psi(V)$ with $\tilde b(\cdot,1) = b$.  Indeed, take a continuous function $\chi:V\to I$ with compact support in $W$, such that $\chi=1$ on a smaller neighbourhood of $C$, and let
\[  \tilde b(s,t)(x) = \left\{ \begin{array}{cl} c\big(s, 1+(t-1)\chi(\pi(x))\big)(x) & \textrm{if $x\in \pi^{-1}(W)$,} \\ b(s)(x) & \textrm{if $x\in \pi^{-1}(V\setminus W)$.} \end{array} \right. \qedhere \]
\end{proof}

Here is the next part of the proof of Theorem \ref{t:main-theorem}(a).

\begin{theorem}   \label{t:local-weak-eq}
The inclusion $\Phi\hookrightarrow \Psi$ is a local weak homotopy equivalence.
\end{theorem}

\begin{proof}
We begin with a self-contained proof, assuming that $X$ is smooth.  Afterwards we consider the more difficult case in which $X$ may be singular.  Let $R$ be a Kempf-Ness set in $X$ as in Section \ref{sec:Stein.compacts}.  Take a point $q\in Q$ and a point $x\in R$ in the closed $G$-orbit in $\pi^{-1}(q)$.  Let $H=G_x$ and $L=K_x$, so $H=L^\C$.  We apply slice theory to the $G$-space $X$ and the $K$-space $R$ at $x$ and obtain arbitrarily small neighbourhoods $U$ of $q$ such that the following hold.
\begin{itemize}
\item  $\pi^{-1}(U)$ is $G$-biholomorphic to $G\times^H S$.  Since $X$ is smooth, the slice $S$ can be chosen to be an $H$-invariant star-shaped neighbourhood of the origin in the $H$-module $T_x X/T_x(Gx)$ \cite[Section 5.5]{Heinzner1991}, so $S$ is holomorphically $H$-contractible, meaning that the identity map of $S$ can be joined to the constant map with value $x$ by a continuous path of holomorphic $H$-maps $S\to S$.  (When $X$ is not smooth, we do not know whether the slice $S$ can be chosen to be holomorphically $H$-contractible.)
\item  $\pi^{-1}(U)\cap R$ is real-analytically $K$-isomorphic to $K\times^L T$.  The slice $T$ is a real-analytic $L$-variety, so it possesses an $L$-equivariant triangulation \cite{Illman2000} and is therefore topologically locally $L$-contractible at the $L$-fixed point $x$, meaning that (after shrinking $U$), the identity map of $T$ can be joined to the constant map with value $x$ by a continuous path of continuous $L$-maps $T\to T$.
\end{itemize}
By adjunction, the restriction maps 
\[ \Phi(U)=\O^G(\pi^{-1}(U), Y) \to \O^H(S, Y) \] 
and 
\[ \mathscr C^K(\pi^{-1}(U)\cap R, Y) \to \mathscr C^L(T, Y) \]
are homeomorphisms.  Moreover, the space of constant maps to $Y^H$ is a deformation retract of each of the spaces $\O^H(S, Y)$ and $\mathscr C^L(T, Y)$.  Finally, since $R$ is a deformation $K$-retract of $X$, $\mathscr C^K(\pi^{-1}(U)\cap R, Y)$ is a deformation retract of $\Psi(U) = \mathscr C^K(\pi^{-1}(U), Y)$.  This shows that $\Phi(U)$ and $\Psi(U)$ both deformation-retract, each in its own way, onto the common subspace of constant maps to $Y^H$.  It follows that the inclusion $\Phi(U) \hookrightarrow \Psi(U)$ is a homotopy equivalence and the proof is complete by Lemma \ref{l:homotopy-equiv}.

In general, when $X$ is not necessarily smooth, we let $\Psi_0(U)$, for $U\subset Q$ open, be the space of continuous $K$-maps $\pi^{-1}(U)\cap R \to Y$ and note that the restriction map $\Psi(U) \to \Psi_0(U)$ is a homotopy equivalence.  Hence, it suffices to show that the morphism $\Phi\hookrightarrow\Psi\to\Psi_0$ is a local weak homotopy equivalence.  With the inclusion $\partial\B_n \hookrightarrow \B_n$ in the definition of a local weak homotopy equivalence replaced by the inclusion of a point in an arbitrary compact Hausdorff space, this is a special case of \cite[Proposition 3.1]{KLS2018}.  The proof of the Proposition is easily adapted to the former inclusion.
\end{proof}

In the proof for the smooth case, we contracted in the source.  We don't know how to do this in the singular case.  In the more intricate argument following the proof of \cite[Proposition 3.1]{KLS2018}, we contract in the target.

\section{Equivariant parametric homotopy approximation} 
\label{sec:approximation}

\noindent
As before, we let $G$ be a reductive complex Lie group, $K$ be a maximal compact subgroup of $G$, and $X$ be a Stein $G$-space.  In this section, only the $K$-action on $X$ is relevant.

\begin{definition}   \label{def:EPHAP}
Let $Y$ be a $K$-manifold with a metric $d$ giving its topology.  Let $L$ be a compact $K$-stable $\O(X)$-convex set in $X$ and let $U\Subset X$ be a $K$-stable Stein neighbourhood of $L$. Suppose that $P$ is a finite polyhedron that deformation-retracts to a point (equivalently, $P$ is contractible) and $P_0$ is a subpolyhedron of $P$.  Set $Q=P\times I$ and $Q_0=(P\times{0})\cup (P_0\times I)$.  Let $f: Q\times X\to Y$ be a continuous $K$-equivariant map such that:
\begin{enumerate}
\item[(i)] for every $q=(p,t)\in Q$, $f_q=f(q,\cdot): X\to Y$ is holomorphic on $U$,
\item[(ii)] for every $q\in Q_0$, $f_q$ is holomorphic on $X$.
\end{enumerate}
We say that $Y$ has the \emph{equivariant parametric homotopy approximation property,} abbreviated EPHAP, if for any $f$ as above and 
$\epsilon>0$, there is a continuous $K$-equivariant map $\tilde f : Q\times X\to Y$ such that for each $q\in Q$, $\tilde f(q,\cdot): X\to Y$ is  holomorphic and:
\begin{enumerate}
\item $\tilde f_q=f_q$ for $q\in Q_0$,
\item  $\underset{x\in L,\ q\in Q}{\sup} d(\tilde f_q(x),f_q(x))<\epsilon$.
\end{enumerate}
\end{definition}

The following theorem is one of the main results of this paper.

\begin{theorem}  \label{thm:K-ell-implies-EPHAP}
Every $K$-elliptic manifold satisfies EPHAP. 
\end{theorem}

\begin{remark}
Suppose that $Y$ is a $G$-manifold and $Y$ satisfies EPHAP (with respect to $K$).  Then the maps $\tilde f_q$ are automatically $G$-equivariant.
\end{remark}

\begin{remark}\label{rem:not-defined-everywhere}
The reader will notice that the conditions (1) and (2) in Definition \ref{def:EPHAP}  do not involve the values of $f_q(x)$ for $x\not\in U$ and $q\not\in Q_0$. In fact, our proof of the theorem shows that one can obtain a family $\tilde f_q$ as required if $f_q$ is a continuous family of $K$-equivariant holomorphic functions from $U$ to $Y$ which extend to be holomorphic $K$-equivariant functions on $X$ for $q\in Q_0$. This is an equivariant version of \cite[Theorem 4.2 and following Remarks]{FP2000}.
\end{remark}

We begin with preliminaries.  Let $E$ be a holomorphic $K$-vector bundle over the Stein $K$-space $X$ with $K$-invariant norm $\vert\cdot\vert$.  The following theorem is the equivariant version of \cite[Theorem 2.8.4]{Forstneric2017}.  Only in this result we allow $P$ to be an arbitrary compact Hausdorff space and $P_0$ to be a closed subspace.

\begin{theorem}[Equivariant Cartan-Oka-Weil theorem with parameters]  \label{thm:C-O-W}
Let $L$ be a compact $\O(X)$-convex $K$-stable subset of $X$ and $X'$ be a $K$-stable subvariety of $X$. Let $\pi_X: P\times X\to X$ be the projection.  Let $f$ be a continuous $K$-equivariant section of $\pi_X^*E$ over $P\times X$ with the following properties.
\begin{enumerate} 
\item[(i)] There is a $K$-neighbourhood $V\Subset X$ of $L$ such that for every $p\in P$, $f(p,\cdot)$ is holomorphic on $V$ and on $X'$.
\item[(ii)] $f(p,\cdot)$ is holomorphic on $X$ for every $p\in P_0$.
\end{enumerate}
Then for every $\epsilon>0$, there is a continuous $K$-equivariant section $F$ of $\pi_X^*E$ such that:
\begin{enumerate}
\item  $F(p,\cdot)$ is holomorphic on $X$ for all $p\in P$,
\item  $\vert F-f\vert<\epsilon$ on $P\times L$,
\item  $F=f$ on $(P_0\times X)\cup (P\times X')$.
\end{enumerate}
\end{theorem}

\begin{proof}
By \cite[Theorem 2.8.4]{Forstneric2017}, there is $F$ with the stated properties, but it might not be $K$-equivariant.  Averaging over $K$ gives an equivariant solution (see \cite[Chapter 3, Theorem 5.1]{BtD1985}).
\end{proof}

\begin{remark}   \label{allow.F.to.change}
If $P$ deformation-retracts to a point, then by Lemma \ref{lem:deform}, the theorem also holds if $\pi_X^*E$ is replaced by a continuous family of holomorphic $K$-vector bundles over $P\times X$.
\end{remark}

Composed sprays will help us give an understandable proof of Theorem \ref{thm:K-ell-implies-EPHAP}.

\begin{definition}
Let $s_1: E_1\to Y$ and $s_2: E_2\to Y$ be $K$-equivariant dominating holomorphic sprays, where $E_1$ and $E_2$ are holomorphic $K$-vector bundles over $Y$ with projections $\pi_1$ and $\pi_2$, respectively.
\begin{enumerate}
\item The composed spray $s_1*s_2: E_1*E_2\to Y$ is defined by
\[  E_1*E_2=\{(e_1,e_2)\in E_1\times E_2: s_1(e_1)=\pi_2(e_2)\}, \]
\[  \pi_1*\pi_2(e_1,e_2)=\pi_1(e_1),\quad s_1*s_2(e_1,e_2)=s_2(e_2).  \]
\item Let $s: E\to Y$ be a $K$-equivariant dominating spray, where $\pi: E\to Y$ is the projection.  For $k\geq 2$, the $k$-th iterated spray map $s^{(k)}: E^{(k)}\to Y$ is defined by
\begin{align*}
E^{(k)}=\{e=(e_1,\dots,e_k): \ &e_j\in E\text{ for }j=1, \dots,k, \\ &s(e_j)=\pi(e_{j+1})\text{ for }j=1,\dots,k-1\},
\end{align*}
\[ \pi^{(k)}(e)=\pi(e_1), \quad s^{(k)}(e)=s(e_k). \]
\end{enumerate}
\end{definition}
Note that $E_1*E_2$ is the pullback of $E_2$ by the spray map $s_1: E_1\to Y$, and similarly for $E^{(k)}$.  Since all the maps involved are $K$-equivariant, $E^{(k)}$ has a holomorphic $K$-action and $s^{(k)}: E^{(k)}\to Y$ is $K$-equivariant.  The bundles $E^{(k)}$ are not naturally $K$-vector bundles over $Y$, but they do have a natural zero section $\Theta=\{(0,\dots,0)\}\subset  E^{(k)}$.  Since $s$ is dominating, so is the differential of $s^{(k)}$ along the fibres of $\pi^{(k)}$ along $\Theta$.

Let $X$ be a Stein $K$-space and let $Z=X\times Y$ with projection $\pi_Y$ to $Y$. 

\begin{proposition}   \label{prop:flatten}
Let $\Omega\subset Z$ be a $K$-stable Stein subset which is either open or a subvariety.  Let $k\geq 2$.  Then there is a $K$-equivariant fibre-preserving biholomorphic map
\[  \Theta: \pi_Y^*E^{(k)}\vert_\Omega\to \oplus^k \pi_Y^*E\vert_\Omega \]
which preserves the zero sections and whose differential at the zero section is the identity.
\end{proposition}

\begin{proof}
This is included in \cite[Lemma 6.3.7]{Forstneric2017} in the non-equivariant case.  The proof is via maps that are automatically $K$-equivariant in our situation.  One also needs the fact that holomorphic $K$-vector bundles over $\Omega$ which are topologically $K$-isomorphic are also holomorphically $K$-isomorphic. This is proved in \cite[\S 11]{HK1995}.
\end{proof}

Let  $P_0$, $P$ and  $f: Q\times X\to Y$, etc.\  be as in Definition \ref{def:EPHAP}. We may assume that $P\subset\R^n\subset\C^n$ for some $n$ so that $Q\subset\R^n\times\R\subset\C^n\times\C$.

\begin{remark}   \label{rem:P_0?} 
It follows from our assumptions that $P_0$ has a neighbourhood $U$ and a deformation retraction $\rho: U\times I\to U$ onto $P_0$.  Using this one can find arbitrarily small neighbourhoods $P_0'$ of $P_0$ and continuous maps $\tau: P\to P$ such that $\tau$ is the identity on a neighbourhood of $U^c$ and $\tau\vert_{P_0'}$ is a retraction onto $P_0$.  If $P_0'$ is sufficiently small, then $f_{(p,t)}(x)$ is arbitrarily close to $f_{(\tau(p),t)}(x)$ on $P\times I\times L$. Thus we may  assume that $f_{q}$ is holomorphic on $X$ for $q\in Q_0':=(P\times{0})\cup (P_0'\times I)$.  In the proofs that follow we may shrink $P_0'$. 
\end{remark}

Let $Z= \C^n\times\C\times X\times Y$ and $\pi_Y: Z\to Y$ the projection. Let $s:  E\to Y$ be the dominating $K$-equivariant spray on the $K$-elliptic manifold $Y$. Let $F:=\pi_Y^*E$ and let $\Phi: F\to Z$ send $(p,t,x,e_y)$  to $(p,t,x,s(e_y))$.

\begin{proposition}   \label{prop:xi.exists}
Let $M:=\gamma_f(Q\times L)$ and $\Omega$ be a Stein $K$-neighbourhood of $M$ in $Z$.  Let $U'\Subset U$ be a Stein $K$-neighbourhood of $L$ such that $\gamma_f(Q\times \overline{U'}) \subset \Omega$.  Let $V\Subset X$ where $V$ is open, $K$-stable, Stein and contains $U'$.  Then any $t_0\in I$ admits a neighbourhood $I_0\subset I$ and a continuous family $\xi_{(p,t)}$ of $K$-equivariant holomorphic sections of  $F\vert_{\gamma_{f_{(p,t_0)}}(U')}$ such that
\begin{equation}\tag{$*$}
\Phi(\xi_{(p,t)})=\gamma_{f(p,t)}  \text{ over $P\times I_0\times U'$, }
\end{equation}
where $\xi_{(p,t_0)}= 0$.  Moreover, shrinking $P_0'$, we can arrange that for $p\in P_0'$, the sections $\xi_{(p,t)}$ extend to be holomorphic on $V$ such that $(*)$ holds with $U'$ replaced by $V$.
\end{proposition}

\begin{proof}
Consider the restriction $\tilde F$ of $F$ to $\Omega$.  Then we have a splitting $\tilde F=\tilde F'\oplus\tilde F''$. For $p\in P$ let $S_p =\gamma_f(p,t_0,U')$.  Then $\Phi$ gives a biholomorphism of a Stein $K$-neighbourhood of the zero section of $\tilde F'\vert_{S_p}$ and a Stein $K$-neighbourhood $\Theta_p$ of $\gamma_f(p,t_0,U')$.  Since $P$ is compact, there is a neighbourhood $I_0\subset I$ of $t_0$ such that $\gamma_ f(p,t,U')\subset\Theta_p$ for $p\in P$ and $t\in I_0$.  Applying the inverse of $\Phi$ we obtain the $\xi_{(p,t)}$ satisfying  $(*)$.  Let $P_0''\Subset P_0'$ be a neighbourhood of $P_0$.  Then we obtain $(*)$  for $U'$ replaced by $V$ and $P$ replaced by $\overline{P_0''}$.  Now using a cutoff function on $P$, we can combine the sections $\xi_{(p,t)}$ of $\tilde F'$ over $\gamma_f(P\times I_0\times U')$ and the sections of $\tilde F'$ over $\gamma_f(\overline{P_0}'' \times I_0\times V)$ to obtain our desired result with $P_0'$ replaced by a neighbourhood of $P_0$ with closure in $P_0''$.
\end{proof}

For $k\geq 1$, let $\Phi^{(k)}:\pi_Y^*E^{(k)}\to Z$ send $(p,t,x,e^{(k)}_y)$ to $(p,t,x,s^{(k)}(e^{(k)}))$. Let $S_p=\gamma_f(p,0,U')$.

\begin{corollary}   \label{cor:main}
There is $k\geq 1$ and a continuous family of $K$-equivariant holomorphic sections $\xi_{(p,t)}$ of $\pi_Y^*E^{(k)}\vert_{S_p}$ such that:
\begin{enumerate}
\item $\xi_{(p,0)}$ is the zero section for each $p\in P$,
\item $\xi_{(p,t)}$ extends to a holomorphic section of $\pi_Y^*E^{(k)}\vert_{\gamma_{f(p,0,V)}}$ for $p$ in a neighbourhood $P_0'$ of $P_0$,
\item $\Phi^{(k)}(\xi_{(p,t)}(\gamma_f(p,0,x)))=\gamma_f(p,t,x)$ for $x\in U'$,  $p\in P$, 
\item for $p\in P_0'$, the above holds for $x\in V$.
\end{enumerate}
\end{corollary}

\begin{proof}
By compactness of $I$, there are numbers $0=t_0<t_1<\dots<t_k=1$ such that for $j=0$, $1, \dots,k-1$, there is a homotopy $\xi^j_{(p,t)}$ of holomorphic sections of $F\vert_{\gamma_f(p,t_j,U')}$ such that
\[  \Phi(\xi^j_{(p,t)}(\gamma_f(p,t_j,x)) )=\gamma_f(p,t,x),\quad t_j\leq t\leq t_{j+1},\ x\in U',\ p\in P.  \]
In particular, $\Phi(\xi^j_{(p,t_{j+1})}(\gamma_f(p,t_j,x)))=\gamma_f(p,t_{j+1},x)$, $j=0,\ 1,\dots,k-1$.  It follows that we can combine the $\xi^j_{(p,t)}$ into a holomorphic section of $\pi_Y^*E^{(k)}\vert_{S_p}$ such that (3) holds.
 For $p\in P_0'$, the sections extend to sections of $\pi_Y^*E^{(k)}$ over $P_0'\times I\times V$ and (4) holds.
 \end{proof}

\begin{proof}[Proof of Theorem \ref{thm:K-ell-implies-EPHAP}]
We may assume that $f_{(p,t)}(x)$ is holomorphic for $p$ in a neighbourhood $P_0'$ of $P_0$.  Let $Q_0' =(P\times 0)\cup (P_0'\times I)$.  Using Proposition \ref{prop:HK}, we find an exhaustion of $X$ by $K$-stable Runge Stein subsets $ W_1\Subset W_2\Subset\cdots\Subset X$ such that each $L_m:=\overline{W_m}$ is $\O(X)$-convex.  We may assume that $U\Subset W_1$.  We show that for any $\epsilon>0$, there is a continuous family of $K$-equivariant holomorphic maps $f^{(1)}: Q\times W_1\to Y$ such that (perhaps shrinking $P_0'$),
\begin{itemize}
\item  $f^{(1)}_q=f_q$ on $Q_0'\times W_1$,
\item $d(f^{(1)}_q(x),f_q(x))<\epsilon/2$ on $Q\times L$.
\end{itemize}
By the same argument, there is a continuous family of $K$-equivariant holomorphic maps $f^{(m)}: Q\times W_m\to Y$, $m\geq 2$, such that
\begin{itemize}
\item  $f^{(m)}_q=f_q$ on $Q_0'\times W_m$,
\item $d(f^{(m)}_q(x),f^{(m-1)}_q(x))<\epsilon/2^{m}$ on $Q\times L_{m-1}$.
\end{itemize}
As $m\to\infty$, the $f^{(m)}$ converge to a continuous $K$-equivariant map $\tilde f$ satisfying the theorem. Thus it is enough to show the existence of $f^{(1)}$. 

By Corollary \ref{cor:main}, there is $k\geq 1$ and a continuous family of  holomorphic sections $\xi_{(p,t)}$ of $\pi_Y^*E^{(k)}\vert_{S_p}$.  For $q\in Q_0'$, $\xi_{(p,t)}$ extends to a $K$-equivariant holomorphic section of $\pi_Y^*E^{(k)}\vert_{\gamma_f(p,0,W_1)}$.  We have $\Phi^{(k)}(\xi_{(p,t)}(\gamma_f(p,0,x)))=\gamma_f(p,t,x)$ for $(p,t,x)$ in $Q\times U'$ and in $Q_0'\times W_1$.   Now the bundle $E^{(k)}\vert_{\gamma_f(p,0,W_1)}$ has the structure of a $K$-equivariant holomorphic vector bundle over $\gamma_f(p,0,W_1)$.  Using this structure and a cutoff function and shrinking $P_0'$, we can extend the $\xi_q$ to be continuous sections defined over $W_1$, unchanged on a neighbourhood of $L$ and unchanged for $q\in Q_0'$.
Using Remark \ref{allow.F.to.change} and Theorem \ref{thm:C-O-W}, we can then find a continuous family of holomorphic sections $\tilde\xi_{(p,t)}$ of $\pi_Y^*E^{(k)}$ over $\gamma_{f_{(p,0)}}(W_1)$ such that $\gamma_{f^{(1)}_q}(x)= \Phi^{(k)}(\tilde\xi_{(p,t)}(\gamma_f(p,0,x)))$, $x\in W_1$, has the required properties.
\end{proof}

\section{Equivariant nonlinear splitting lemma}
\label{sec:splitting}

\noindent
In this section, we generalise to an equivariant setting a version of the nonlinear splitting lemma \cite[Proposition 5.8.4]{Forstneric2017} that first appeared in  \cite{Forstneric2007} and \cite{DF2007}. Our arguments follow closely those of \cite{FP2000} and \cite{FP2001}.  As before, we let $G$ be a reductive complex Lie group, $K$ be a maximal compact subgroup of $G$, and $X$ be a Stein $G$-space.  Here, as in the previous section, only the $K$-action on $X$ is relevant.  Let $\pi:X\to X\git G = X\git K$ be the categorical quotient map.  We first recall the following definitions (see \cite[Definitions 5.7.1 and 6.7.1]{Forstneric2017}).
 
\begin{definition}   \label{C-pairs}   
A compact subset $A$ of a complex space $S$ is a Stein compact if it admits a basis of Stein neighbourhoods in $S$.  Let $A$, $B$ be compact sets in $S$.  The pair $(A,B)$ is called a Cartan pair if
\begin{enumerate}
\item $A$, $B$, $C=A\cap B$ and $D=A\cup B$ are Stein compact subsets of $S$,
\item $\overline{A\setminus B}\cap \overline{B\setminus A}=\varnothing$.
\end{enumerate}
The pair $(A,B)$ is called a $\CC$-pair if, in addition,
\begin{enumerate}
\addtocounter{enumi}{2}
\item $C$ is $\O(B)$-convex.
\end{enumerate}
\end{definition}

\begin{remark}
An $\O(X)$-convex compact subset  in a Stein space is a Stein compact.
\end{remark}

Let $\phi: X\to [0,\infty)$ and the associated Kempf-Ness set $R$ be as in Section \ref{sec:Stein.compacts}.  

\begin{lemma}  \label{lem:lifting-Cartan-pairs}
Let $(A_0,B_0)$ be a Cartan pair in $X \git K$.  We construct a Cartan pair $(A, B)$ of $K$-invariant subsets of $X$ such that $\pi^{-1}(A_0)\cap R\subset A$ and $\pi^{-1}(B_0)\cap R\subset B$, so $\pi(A)=A_0$ and $\pi(B)=B_0$.  If $(A_0,B_0)$ is a $\CC$-pair, then $(A,B)$ is a $\CC$-pair.
\end{lemma}

\begin{proof}
Let $D_0=A_0\cup B_0$.  Let $r=\sup\phi(x)$ for $x\in \pi\inv(D_0)\cap R$.  It follows from Proposition \ref{prop:HK} that $\tilde D:=\phi\inv([0,r])$ is $K$-invariant  and $\O(X)$-convex.  Since intersections of $\O (X)$-convex subsets are $\O (X)$-convex, $A:=\pi\inv(A_0)\cap \tilde D$ and $B:=\pi\inv(B_0) \cap\tilde D$ form a $K$-invariant Cartan pair in $X$ satisfying the lemma. If $A_0\cap B_0$  is $\O(B_0)$-convex, then  $\pi\inv(A_0\cap B_0)\cap \tilde D$ is $\O(\pi\inv(B_0)\cap\tilde D)$-convex, that is, $A\cap B$ is $\O(B)$-convex. 
\end{proof}
 
Let $Z$ be a Stein $K$-manifold and $D$ a domain in $Z$.  We denote by $H^\infty(D)$ the Banach algebra of bounded holomorphic  functions on $D$. 

\begin{lemma}  \label{lem:cartan.open.sets}
Let $(A,B)$ be a Cartan pair of $K$-stable sets in $Z$ where $C=A\cap B\neq \emptyset$.  Then there are arbitrarily small $K$-stable Stein neighbourhoods $A'$ and $B'$ of $A$ and $B$, respectively, satisfying:
\begin{enumerate}
\item  $\Omega=A'\cup B'$ has smooth boundary and is strongly pseudoconvex,
\item  $\overline{A'\backslash B'} \cap \overline{B'\backslash A'} =\emptyset$.
\end{enumerate}
\end{lemma}

\begin{proof}
Let $\widetilde A\supset A$ and $\widetilde B\supset B$ be Stein neighbourhoods of $A$ and $B$. By the proof of \cite[Lemma 2.4]{FP2000}   there are smaller  Stein open sets $A\subset A_0\subset\widetilde A$ and $B\subset B_0\subset \widetilde B$ satisfying (2). By Lemma \ref{lem:intersection.Stein} we may assume that $A_0$ and $B_0$ are $K$-invariant and we may choose a $K$-invariant Stein neighbourhood $\Omega'$ of $A\cup B$ in $A_0\cup B_0$. Choose a smooth strictly plurisubharmonic $K$-invariant exhaustion $\phi: \Omega'\to\R^+$. Then for $d>0$ sufficiently large  $\Omega=\{\phi<d\}$ contains $A\cup B$. Choose $d$ so that it is also a regular value of $\phi$. Then $\partial\Omega$ is smooth and $\Omega$ is  strongly pseudoconvex and $K$-invariant. The  argument of \cite[Lemma 2.4]{FP2000} now shows that $A':=\tilde A\cap\Omega$ and $B':=\tilde B\cap \Omega$ have union $\Omega$ and the lemma is established. 
 \end{proof}

\begin{lemma}  \label{lem:linear.split} 
Let $A, B \Subset Z$ be relatively compact $K$-domains satisfying the following:
\begin{enumerate}
\item  $\Omega=A\cup B$ is a smooth strongly pseudoconvex domain in $Z$ and
\item  $\overline{A\backslash B} \cap \overline{B\backslash A} =\emptyset$.
\end{enumerate}
Let $C=A\cap B$. Then there exist bounded linear operators
${\calA}: H^\infty(C)^K \to H^\infty(A)^K$ and 
${\calB}: H^\infty(C)^K \to H^\infty(B)^K$ 
such that $c={\calA}(c) - {\calB}(c)$ on $C$ for each
$c\in H^\infty(C)^K$.
\end{lemma}

\begin{proof}
The non-equivariant version of the above is \cite[Lemma 3.2]{FP2001} and we get our desired result by averaging over $K$.
\end{proof}

\begin{remark}  \label{rem:bundles} 
Let $V$ be a $K$-module and let $H^\infty(A,V)$ denote the space of bounded holomorphic functions on $A$ with values in $V$.  Then Lemma \ref{lem:linear.split} applies with $H^\infty(A)^K$ replaced by $H^\infty(A,V)^K$ and similar replacements for $H^\infty(B)^K$ and $H^\infty(C)^K$. 
\end{remark}

Let $V_r$  be the ball of radius $r$ in a unitary $K$-module $V$. Let $D\Subset Z$ be a relatively compact $K$-domain. Let $H^\infty(D\times V_r,V)$ denote the space of bounded holomorphic maps of $D\times V_r$ to $V$ and $H^\infty(D\times V_r,V)^K$ the subset of $K$-equivariant ones.

\begin{proposition}  \label{prop:main-splitting}
Let $(A,B)$ be a $K$-invariant Cartan pair in $Z$ and let $C=A\cap B$.  Let $r>0$   and let $\tilde C\Subset Z$ be a $K$-stable neighbourhood of $C$.  Let $\psi_0: \widetilde{C}\times V_r\to V$ be projection on the second factor.  Then there are $K$-neighbourhoods $A'\supset A$ and $B'\supset B$ with $C'=A'\cap B'\Subset   \widetilde C$ (where $A'$ and $B'$ satisfy conditions (1)--(2) of Lemma \ref{lem:linear.split}) and an open neighbourhood ${\calW}$ of $\psi_0$ in the Banach space $H^\infty(\widetilde C\times V_r,V)^K$ and smooth Banach
space operators
${\calA'}: {\calW}\to H^\infty(A',V)^K$,
${\calB'}: {\calW}\to H^\infty(B',V)^K$,
with ${\calA'}(\psi_0)=0$ and ${\calB'}(\psi_0)=0$. 
For each $\psi\in {\calW}$ the  bounded holomorphic maps
$\alpha={\calA'}(\psi) : A'\to V$ and
$\beta={\calB'}(\psi) :  B'\to V$ satisfy
$$
   \psi\bigl(x, \alpha(x) \bigr)= \beta(x) \quad (x\in C').                                                    
$$
If $\psi \in {\calW}$ satisfies $\psi(x,0)=0$ for
$x\in \widetilde C$, then ${\calA'}(\psi)=0$ and ${\calB'}(\psi)=0$.
\end{proposition}

\begin{proof}
The non-equivariant case is  \cite[Proposition 4.1]{FP2001}.  However, since the maps $\calA'$ and $\calB'$ are not linear, we cannot use an averaging  argument to get the above.  Let  ${\calA}: H^\infty(C',V)^K \to H^\infty(A',V)^K$ and ${\calB}: H^\infty(C',V)^K \to H^\infty(B',V)^K$ be as in Lemma \ref{lem:linear.split} (see Remark \ref{rem:bundles}).  Consider the operator
\begin{align*}  & \Phi: H^\infty(\widetilde{C}\times V_r,V)^K \times H^\infty(C',V)^K 
   \to H^\infty (C',V)^K, \\
    & \Phi(\psi,c)(x) = \psi\bigl(x,{\calA}c(x) \bigr) - {\calB}c(x),
    \qquad x\in C'.  
\end{align*}
By the argument in \cite[Proposition 5.2]{FP2000}, $\Phi$ is a smooth Banach space operator in a neighbourhood of $(\psi_0,0)$. Note that 
$$
\Phi(\psi_0,c)(x)=\psi_0(x,(\calA c)(x))-\calB(c)(x)=c(x), \quad x\in C'.
$$
Thus the partial derivative $\partial_c\Phi(\psi,c)$ at $(\psi_0,0)$ is the identity.  By the implicit function theorem, the equation $\Phi(\psi,c)=0$ has a unique solution $c=\calC(\psi)$ in a neighbourhood of $\psi_0$ where $\calC$ is a smooth operator.  Then the operators $\calA'=\calA\circ\calC$ and $\calB'=\calB\circ\calC$ establish the proposition.
\end{proof}

Let $X$ be a Stein $K$-space and let  $r_0>r>0$.  Let $\lVert \cdot \rVert$ be the $K$-invariant norm on $V$.  Let $U\Subset X$ be open and $K$-stable and let $\gamma: U\times V_{r_0}\to V$ be holomorphic.  Let $\id$ be the projection $U\times V_{r_0}\to V$. Write $\gamma=\id+\gamma'$ and let 
$$
\lVert\gamma'\rVert=\sup_{x\in U,v\in V_{r_0}}\lVert\gamma'(x,v)\rVert.
$$

\begin{theorem}\label{thm:main-splitting}
Let $(A,B)$ be a $K$-invariant Cartan pair in $X$ and let $C=A\cap B$, $D=A\cup B$. Let $U$, $r$, and $r_0$ be as above.  Then there exist arbitrarily small open $K$-neighbourhoods $A' \supset A$, $B' \supset B$ with $C'=A'\cap B' \subset U$ and a number $\delta>0$ satisfying the following.  For every  holomorphic $K$-map 
$$
\gamma : U\times V_{r_0}\to V\text{ with }\lVert\gamma'\rVert<\delta 
$$
there exist holomorphic $K$-maps $\alpha_\gamma: A'\to V$ and $\beta_\gamma : B'\to V$, depending continuously on $\gamma$, with $\alpha_\id =0$ and $\beta_\id=0$, satisfying
$$
\gamma(x,\alpha_\gamma(x))=\beta_\gamma(x),\quad x\in C'.
$$
\end{theorem}

\begin{proof}  
We follow part of the argument of \cite[Proposition 5.8.4]{Forstneric2017}.  Since $D = A \cup B$ is a $K$-invariant Stein compact in $X$, there is an open $K$-invariant Stein neighbourhood $U_0 \Subset X$ of $D$ which we may assume contains $U$.  By Proposition \ref{prop:equivariant-embedding}(b), $U_0$ embeds $K$-equivariantly as an open subset of a $G$-subvariety of a $G$-module $W \simeq \C^m$.  Replace $X$ by $U_0$ and identify $X$ with its image in $W$ so that $D$ becomes a $K$-invariant Stein compact in $W$.  We may assume that $U\supset C$ equals $X \cap \tilde{U}$, where $\tilde U\Subset W$ is a $K$-stable Stein domain.  Choose a $K$-stable Stein domain $\tilde C\Subset\tilde U$ such that $C\subset \tilde C\cap X\Subset U$.  Applying  \cite[Lemma 2.8.3]{Forstneric2017} with $\Omega=\tilde  C\times V_r$ and $\Omega'=\tilde U\times V_{r_0}$ with $X'=X\times V_{r_0}$, we obtain a bounded linear extension operator 
$$
S: H^\infty(U\times V_{r_0},V)\to H^\infty(\tilde C\times V_r,V).
$$ 
By averaging over $K$ we obtain a bounded linear extension operator $\tilde S$ for the $K$-invariant elements.  Set
$$
\tilde\gamma(x,v)=v+\tilde S (\gamma')(x,v),\quad x\in\tilde C,\ v\in V_r.
$$
By properties of $S$ and $\tilde S$, we have $\gamma=\tilde\gamma\vert_{(\tilde C\cap X)\times V_r}$.  Since $\tilde S$ is continuous, we have $\Vert \tilde S(\gamma')\rVert<M\delta$, where $M$ is the norm of $\tilde S$.

Let $\tilde D\Subset W$ be a $K$-stable open set in which $X$ is closed. Choose $K$-stable domains $\tilde A$, $\tilde B$ in $W$, both contained in $\tilde D$,  such that $\tilde A\cap\tilde B\subset \tilde C$ and $\tilde A$ and $\tilde B$ satisfy properties (1) and (2) of Lemma \ref{lem:cartan.open.sets}.  Then Proposition \ref{prop:main-splitting} gives us equivariant maps
$$
\alpha_{\tilde \gamma}: \tilde A\to V \text { and } \beta_{\tilde \gamma}: \tilde B\to V \text{ with } \tilde\gamma(x,\alpha_{\tilde\gamma}(x))=\beta_{\tilde\gamma}(x)
$$
as long as $M\delta$ is sufficiently small.  Let $\alpha_\gamma$ be the restriction of $\alpha_{\tilde\gamma}$ to $A':=\tilde A\cap X$ and similarly define $B'$ and $\beta_\gamma$.  Then $\gamma(x,\alpha_\gamma(x)) =\beta_\gamma(x)$ for $x\in C'=A'\cap B'$, establishing the theorem.
\end{proof}

\begin{remark}\label{rem:params}
If $\gamma$ depends continuously on parameters in a compact Hausdorff space $P$ with $\gamma=\id$ on a closed subspace $P_0\subset P$, then the theorem gives us $\alpha_\gamma$ and $\beta_\gamma$ depending continuously on $P$ and vanishing on $P_0$.
\end{remark}

\section{Holomorphic weak flexibility} 
\label{sec:hol-weak-flex}

\noindent
In this section, we prove Theorem \ref{t:hol-flex}, the holomorphic weak flexibility of the sheaf $\Phi$ defined in the introduction.  This is the most substantial of the three parts of the proof of Theorem \ref{t:main-theorem}(a).  It follows from an equivariant version of the Heftungslemma \cite[Proposition 6.7.2]{Forstneric2017}. 

As before, we let $G$ be a reductive complex Lie group, $K$ be a maximal compact subgroup of $G$, and $X$ be a Stein $G$-space.  Let $Y$ be $K$-elliptic with dominating spray map $s: E\to Y$.  Other assumptions are as in Definition \ref{def:EPHAP}.  We consider continuous families of $K$-equivariant holomorphic maps $f_p: U\to Y$, $p\in P$, where $U\subset X$ is  $K$-stable and open. We call such maps \emph{holomorphic $P$-families.} 

\begin{definition}
Let $F$ be a holomorphic $K$-vector bundle over $X$. The space of holomorphic sections of $F$ over $X$ is denoted by $\Gamma(X,F)$ (or just $\Gamma(F)$ if $X$ is understood).  Now $K$ acts on $\Gamma(F)$ by 
$$
(k\cdot\xi)(x)=k(\xi(k\inv(x))), \quad k\in K,\  \xi\in\Gamma(F),\ x\in X.
$$
We say that $\xi$ is $K$-finite, and write $\xi\in \Gamma(F)_K$, if $K\cdot\xi$ spans a finite-dimensional subspace of $\Gamma(F)$ (which is automatically $K$-stable).
For a finite-dimensional $K$-module $V$, let $\Gamma(F)_V$ denote the direct sum of the subspaces of $\Gamma(F)_K$ which are $K$-isomorphic to $V$.  Note that $\Gamma(F)_K$ is the direct sum of the spaces $\Gamma(F)_V$ where $V$ runs over  the isomorphism classes of irreducible $K$-modules. Each $\Gamma(F)_V$ is a  module  over $\O(X)^K$.
\end{definition}

From \cite{Harish-Chandra1966} or \cite[Chapter 3, Theorem 5.7]{BtD1985}  we have the following.

\begin{lemma}  \label{lem:K-finite}
Let $F$, etc.\ be as above. Then $\Gamma(F)_K$ is dense in $\Gamma(F)$ in the Fr\'echet topology. This also holds true with $\Gamma(F)$ replaced by any closed $K$-invariant subspace.
\end{lemma}

\begin{corollary}\label{cor:K-finite}
Let $U\Subset X$  be $K$-stable and Stein. Then there is an $N\in\NN$ and  $\xi_1,\dots,\xi_N\in\Gamma(F)_K$  which form a basis for an $N$-dimensional $K$-module $W$  and which generate $\Gamma(U,F)$ over $\O(U)$.  If $X'\subset X$ is a $K$-stable closed subspace with ideal sheaf $I_{X'}$, then $(I_{X'}\Gamma(U,F))_K$  has  finitely many   generators over $\O(U)$ which form  the basis of a $K$-module. 
\end{corollary}

\begin{proposition}   \label{prop:equivariant-6.7.2}
Let $(A,B)$ be a $K$-invariant $\CC$-pair in $X$ and $\tilde A\supset A$ and $\tilde B\supset B$ be $K$-stable open sets.  Assume that $a: \tilde A\times P\to Y$ and $b: \tilde B\times P\to Y$ are holomorphic $P$-families whose restrictions to $\tilde C=\tilde A\cap\tilde B$ are homotopic by a homotopy $c_s$  of holomorphic $P$-families with $c_0=a\vert_{\tilde C}$, $c_1=b\vert_{\tilde C}$. We assume that $c_s$ is constant on $\tilde C\times P_0$.  Then there are $K$-stable open sets $A'$, $B'$, and $C'$ with $A\subset A'\subset\tilde A$, $B\subset B'\subset\tilde B$, and $C\subset C'\subset \tilde C$ such that for any $\epsilon>0$, there are homotopies of holomorphic $P$-families $a_t: A'\times P\to Y$, $b_t: B'\times P\to Y$, $t\in I$, and a homotopy $c_{s,t}$ of the holomorphic $(P\times I)$-family $c:  C'\times P\times I\to Y$ such that
\begin{enumerate}
\item $c_{0,t}=a_t$ and $c_{1,t}=b_t$ over $C'$,
\item $a_t$, $b_t$ and $c_{s,t}$ are independent of $t$ when $p\in P_0$,
\item $c_{s,1}$ is independent of $s$,
\item $d(a_t(x,p),a(x,p))<\epsilon$ for all $x\in A'$, $p\in P$, $t\in I$.
\end{enumerate}
\end{proposition}

\begin{proof}  We may assume that $\tilde A$, $\tilde B$ and $\tilde C$ are relatively compact $K$-invariant  Stein domains.  Note that $c_1$ extends to $b: \tilde B\times P\to Y$.  We now thicken $a$ and $b$ by taking the products of their domains with a small ball.

As before, we may assume that $P\subset \R^n\subset\C^n$. Let $Z=\tilde A\times \C^n\times Y$ with projection  $\pi_Y$ to $Y$.  By Lemma \ref{lem:U.exists}, $M=\gamma_a(A\times P)$ admits a $K$-stable Stein neighbourhood $U\subset Z$. By shrinking $\tilde A$ we may assume that $\overline{\gamma_a(\tilde A\times P)}\subset U$. We have the holomorphic $K$-vector bundle $F:=\pi_Y^*E$ over $U$ with the induced spray map $\sigma:=\pi_Y^*s$.  By Corollary \ref{cor:K-finite}, we may find $\xi_1,\dots,\xi_N\in \Gamma(F)_K$, which form a basis for a $K$-module $W\simeq\C^N$ and generate $F$ at every point of $\overline{\gamma_a(\tilde C\times P)}$.  By construction, $w\in \C^N\mapsto \sum_i w_i\xi_i$ is a $K$-equivariant map. We may assume that $K\to \GL_N(\C)$ has image in $\U_N(\C)$.  Let $W_r$ denote the ball of radius  $r>0$ in  $W$. We have a family $a': \tilde A\times W_r\times P\to Y$,
$$
(x,w,p)\mapsto \sigma(\sum w_i\xi_i(x,p,a(x,p))).
$$
Since $c_0=a$ on $\tilde C$, we may find $0<s_1\leq 1$ such that $(x,p,c_s(x,p))\in U$ for $s\in[0,s_1]$. Choose a continuous function $\chi:I\to I$ which equals $1$ near $0$ and has support in $[0,s_1)$.
Define $c'_s:\tilde C\times W_r\times I\times P\to Y$,
$$
(x,w,s,p)\mapsto \sigma(\sum \chi(s)w_i\xi_i(x,p,c_s(x,p))).
$$
Note that for $x\in\tilde C$,
$$
a'(x,0,p)=a(x,p), \quad c'_s(x, 0,p)=c_s(x,p), \quad a'(x,w,p)=c'_0(x,w,p).
$$
For $s\in[s_1,1]$ and $x\in\tilde C$,
$$
c_s'(x,w,p)=c_s(x,p)\text { so that } c'_1(x,w,p)=b(x,p).
$$  
By Remark \ref{rem:not-defined-everywhere}, using that $C$ is $\O(B)$-convex and perhaps shrinking $\tilde B$, we obtain a homotopy of $K$-equivariant holomorphic $P$-families  $\tilde c'_s: \tilde B\times W_r\times P\to Y$ with $\tilde c'_1(x,w,p)=c_1'(x,w,p)=b(x,p)$, and $\tilde c'_s$ approximates $c'_s$ as closely as desired on  $\tilde C\times W_r\times P$ for $s\in I$.  In particular, when $s=0$, $\tilde c'_0(x,w,p)$ approximates $a'(x,w,p)$ as closely as one wants on $\tilde C\times W_r\times P$.

Replacing our original $c_s'$ by $\tilde c'_s$ we reduce to the case that the homotopy $c_s'$ is defined on all of $\tilde B\times W_r\times P$. We rename it to $b_s'$.  Then  $b_1'=b$ and  $b_0'$ is arbitrarily close to $a'$ on $\tilde C\times  W_r\times P$. We may have had to shrink the open sets around $A$, $B$ and $C$ in our process. Let $a'_s=a'$ be the constant family, $s\in I$.

Shrinking $\tilde C$, making $r$ smaller, using that $a'$ and $b'_0$ are close on $\tilde C\times W_r\times P$ and that the spray map is dominating and equivariant, we can find a $K$-subspace $L\subset W$ such that  
$$
Da'_s(x, l,p)\vert_{L}: L\to T_{a'(x,l,p)}Y \quad\textrm{and}\quad Db'_s(x,l,p)\vert_{L}: L\to T_{b'(x,l,p)}Y
$$
are $K$-isomorphisms for $x\in \tilde C$, $l\in L_r$, $p\in P$, and $s\in I$.

Let $ \Phi_{a'}: \tilde C\times P\times I\times L_r\to \tilde C\times P\times I\times Y$ be given by
$$
(x,p,s,l)\mapsto (x,p,s,a_s'(x, l,p)).
$$
Then $\Phi_{a'}$ is a $K$-equivariant $(P\times I)$-family of local biholomorphic maps in a neighbourhood of $0\in L_r$. 
Similarly define $\Phi_{b'}$, which is also a $K$-equivariant $(P\times I)$-family of local biholomorphisms. Then $\gamma:=\Phi_{b'} \inv\circ \Phi_{a'}$ is a $(P\times I)$-family of local $K$-automorphisms of the $L_r$-bundle $\tilde C\times L_r \times P\times I$ near the identity. For $p\in P_0$,  $\gamma$ is the identity.  As in Theorem \ref{thm:main-splitting}, we write $\gamma=\id+\gamma'$ where $\gamma'$ is close to zero. Let $\gamma_t=\id+t\gamma'$, $0\leq t\leq 1$.  Then $\gamma$ depends upon parameters $(p,s,t)\in P\times I^2$.  We may choose $\tilde C$ and $r$ and $\tilde c'_s$ above so that $\lVert\gamma'\rVert$ is small enough for Theorem \ref {thm:main-splitting} to apply.  Then we obtain $P\times I^2$-families $\alpha_\gamma$ and $\beta_\gamma$ mapping $\tilde C$ to $L$, such that $\gamma(x,\alpha_\gamma(x)) =\beta_\gamma(x)$ over $\tilde C\times P\times I^2$. When $t=0$, $\gamma$ is the identity and $\alpha_\gamma$ and $\beta_\gamma$ vanish. By construction, $\Phi_{a'}\circ\alpha_\gamma=\Phi_{b'}\circ\beta_\gamma$ and hence we may modify $a'$ and $b'$ such that they agree on $\tilde C\times P\times I^2$. We have modified $a$ and $b$ by a family $c_{s,t}$ over a neighbourhood of $C$ with the required properties. The only problem with our construction is (1), since $c_{0,0}$ is only close to $a$ on $C'\times P$, but this can be fixed since $c_{0,0}$ and $a$ are connected by a homotopy over $C'\times P$.
\end{proof}

\begin{theorem}   \label{t:hol-flex}
Every $\CC$-pair $(A_0,B_0)$ in $Q=X\git G$ is weakly flexible for $\Phi$. 
\end{theorem}

\begin{proof}
Let $(A_0,B_0)$ be a $\CC$-pair in $Q$ with neighbourhoods $\tilde A_0$ and $\tilde B_0$ of $A_0$ and $B_0$, respectively.  Let $\tilde C_0\subset\tilde A_0\cap \tilde B_0$ be a neighbourhood of $C_0=A_0\cap B_0$. We are given $P$-families of holomorphic $G$-maps $a:\pi^{-1} (\tilde A_0 )\times P \to Y$ and $b:\pi^{-1} (\tilde B_0) \times P\to Y$ and a homotopy of $P$-families of holomorphic $G$-equivariant maps $c_s: \pi^{-1} (\tilde C_0)\times P \to Y$ between the restrictions of $a$ and $b$. On $\pi\inv(\tilde C_0)\times P_0$ we have $a=b$ and the homotopy $c_s$ is  constant.  By Lemma \ref{lem:lifting-Cartan-pairs}, there is a $K$-invariant $\CC$-pair $(A, B)$ in $X$ such that  $\pi(A)=A_0$ and $\pi(B)= B_0$. By construction, $\pi\inv(A_0)\cap R\subset A$ and similarly for $B$.  Choose a $K$-stable neighbourhood $\tilde A$ of $A$ which is contained in $\pi\inv(\tilde A_0)$, and similarly choose $K$-neighbourhoods $\tilde B$ of $B$ and $\tilde C$ of $C=A\cap B$.
 
We now restrict $a$, $b$ and $c_s$ to the open sets $\tilde A$, $\tilde B$, and $\tilde C$.  By Proposition \ref{prop:equivariant-6.7.2}, replacing $\tilde A$, etc.\ by smaller neighbourhoods $A'$, etc.\ we can find homotopies $a_t$ and $b_t$ connected by a homotopy $c_{s,t}$ satisfying (1)--(4) of the proposition.

The last step is to extend our maps and homotopies to $G\cdot A'$, $G\cdot B'$, and $G\cdot C'$. We are allowed to shrink $A'$, etc.\ to accomplish this. Now by \cite{HK1995}, we may find arbitrarily small $K$-neighbourhoods $U$ of arbitrary compact subsets of the Kempf-Ness set $R$ (see Section \ref{sec:Stein.compacts}) which are orbit convex. For such a neighbourhood $U$, any $K$-equivariant holomorphic map $U\to Y$ extends uniquely to a $G$-equivariant holomorphic map from $G\cdot U$ to $Y$.  By Lemma \ref{lem:lifting-Cartan-pairs}, $A'$ contains a neighbourhood of $A\cap R$ and similarly for $B'$ and $C'$, so we can shrink $A'$, etc.\ and get our desired result.
\end{proof}

This concludes the proof of Theorem \ref{t:main-theorem}(a).

\section{Interpolation and approximation} 
\label{sec:interpolation}

\noindent
In this section, we show how to incorporate interpolation into the proof of Theorem \ref{t:main-theorem}(a) so as to prove Theorem \ref{t:main-theorem}(b) and (d).  As before, we let $G$ be a reductive complex Lie group, $K$ be a maximal compact subgroup of $G$, $X$ be a Stein $G$-space, $\pi:X\to Q=X\git G$ be the categorical quotient, and $Y$ be a $G$-elliptic manifold.  Now we take a $G$-invariant subvariety $X'$ of $X$ and a holomorphic $G$-map $h:X'\to Y$ and redefine the sheaves $\Phi\hookrightarrow \Psi$ on $Q$ by letting $\Phi(U)$, where $U\subset Q$ is open, be the space of holomorphic $G$-maps $f:\pi^{-1}(U)\to Y$ with $f=h$ on $X'\cap\pi^{-1}(U)$, and $\Psi(U)$ be the space of continuous $K$-maps $f:\pi^{-1}(U)\to Y$ with $f=h$ on $X'\cap\pi^{-1}(U)$.  Both spaces are endowed with the compact-open topology.  Theorem \ref{t:main-theorem}(b) states that the inclusion $\Phi(Q)\hookrightarrow\Psi(Q)$ is a weak homotopy equivalence and is proved as follows.

The proof of weak flexibility of $\Psi$ (Proposition \ref{p:top-flex}) holds in the present setting unchanged. To prove weak flexibility of $\Phi$, we need to  incorporate interpolation into Proposition \ref{prop:equivariant-6.7.2} which uses  Theorem \ref{thm:K-ell-implies-EPHAP}, hence we need versions of both with interpolation. Once this is done, the proof of Theorem \ref{t:hol-flex}  goes through unchanged, establishing that $\Phi$ is weakly flexible.  

For Theorem \ref{thm:K-ell-implies-EPHAP} we modify EPHAP (Definition \ref{def:EPHAP}) by adding:
\begin{enumerate}
\item[(iii)] $f_q=h$ on $X'$ for all $q\in Q$,
\item[(3)] $\tilde f_q=h$ on $X'$ for all $q\in Q$.
\end{enumerate}
We have to prove the theorem with these added constraints.  Here are the necessary changes. In the proof of Proposition \ref{prop:xi.exists}, we have the result that $K$-equivariant maps $f_{p,t}$ sufficiently close to some $f_{p,t_0}$, $t_0\in I$, correspond to $K$-equivariant holomorphic sections of a holomorphic $K$-vector bundle $F$ with base the graph of $f_{p,t_0}$.  Since all maps $f_{p,t}$ equal $h$ on $X'\cap U'$, the sections we get all vanish there.  This implies that the constructions in Corollary \ref{cor:main} lead to sections of iterated bundles that all vanish on $X'$. In the proof of Theorem \ref{thm:K-ell-implies-EPHAP} we use Proposition \ref{prop:flatten} to \lq\lq flatten\rq\rq\ the sections of the iterated bundles (which are zero on $X'$) to get sections of vector bundles (which then vanish on $X'$). Then we go back from sections of vector bundles  to maps $X\to Y$. Since all the sections of the vector bundles are zero on $X'$, the functions we construct all equal $h$ on $X'$. Thus Theorem \ref{thm:K-ell-implies-EPHAP} holds with the new conditions (iii) and (3).

To establish Proposition \ref{prop:equivariant-6.7.2} with interpolation, we, of course, change the hypotheses so that $a$, $b$ and $c_s$ are equal to $h$ on the intersection of their domains with $X'\times P$. When we thicken $a$ and $b$ we replace the finite generating set of $K$-finite sections of $F$ over $U$ by a finite set of $K$-finite sections $\xi_i$ of $I_{X'}\Gamma(U,F)$ where $I_{X'}$ is the ideal sheaf of $X'$. Then the family $a':\tilde A\times W_r\times P\to Y$ equals $h$ over $\tilde A\cap X'$ and all the following constructions involve maps which equal $h$ whenever a point of $X'$ is involved. At the end, we are able to modify $a$ and $b$ by a family $c_{s,t}$ over a neighbourhood of $C$ where all maps equal $h$ on the intersection of $X'$ with their domains, giving the  version of Proposition \ref{prop:equivariant-6.7.2} with interpolation.

Next we complete the proof of Theorem \ref{t:main-theorem}(b).  Unlike the proofs in Section \ref{sec:local}, the following proof uses the $G$-ellipticity of $Y$.

\begin{theorem}  \label{t:local-weak-eq-with-interpolation}
The inclusion $\Phi\hookrightarrow \Psi$ is a local weak homotopy equivalence.
\end{theorem}

\begin{proof}
Let $q_0\in Q$ and let $Q'=\pi(X')$. If $q_0\in Q\setminus Q'$, we obtain local weak homotopy equivalence as before.  So assume that $q_0\in Q'$. Let $P=\B_n$ and $P_0=\partial\B_n$.   Let $x_0\in R\cap\pi\inv(q_0)$. Then $Gx_0$ is the closed $G$-orbit in $\pi\inv(q_0)$ and it lies in $X'$. It follows from \cite[4.4 Theorem] {Heinzner1991} that $Kx_0$ is $\O(Gx_0)$-convex, hence $\O(X)$-convex.  Let $Z$ denote $\C^n\times X\times Y$ where $P\subset\R^n\subset\C^n$. As in Section \ref{sec:background}, the graph
\[ M=\gamma_h(P\times Kx_0)=\{(p,kx_0,kh(x_0)) : p\in P,\ k\in K\} \]
is a Stein compact subset of $Z$.  Thus there is a basis of Stein $K$-neighbourhoods $V_k$ of $M$. Let $\pi_Y : Z\to Y$ be projection onto $Y$. The pullback $F=\pi_Y^*E$ has a dominating fibrewise $K$-equivariant spray map $\sigma=\pi_Y^*s: F\to Y$ which for $z=(p,x,y)\in Z$ sends $F_{(p,x,y)}=E_y\to Y$ via $s$.  For any $p\in P$, $f_p(Kx_0)=h(Kx_0)$, so that for any $k\in\NN$ we can choose a Stein $K$-neighbourhood $U_k$ of $Kx_0$ such that $\gamma_{f_p}(P\times U_k)\subset V_k$. There is a subbundle $F'\subset F$ such that $\pi_Y^*s: F'\to Y$ is fibrewise dominating, where the fibre dimension of $F'$ is $\dim Y$. By Theorem \ref{thm:eq-biholom}, if  $V_k$ is sufficiently small, we can lift any $\gamma_{f_p}(P\times U_k)$ to a unique  continuous section $\xi_{p}\in\Gamma(F'\vert_{M})^K$. If $p\in P_0$, then the section is holomorphic. Moreover, every $\xi_p$ is zero on $U_k\cap X'$. The direct sum of $F'$ and another holomorphic $K$-vector bundle over $V_k$ is a product bundle $V_k\times W$, where $W$ is a complex $K$-module. Thus our $f_p$ restricted to $P\times U_k$  correspond to continuous $P$-families of $K$-equivariant sections $\xi_p$ of the product bundle $U_k\times W$, which are holomorphic when $p\in P_0$ and vanish for  $x\in U_k\cap X'$.

Suppose that we can construct a homotopy $\xi_{p,t}$ in our family of maps $\xi_p$ as above, such that when $t=1$, the $\xi_p$ are all holomorphic. Using the $K$-equivariant projection of $U_k\times W$ to $F'\vert_M$ and our spray map, we obtain a homotopy $f_{p,t}$ of our original family $f_p$ of $K$-equivariant maps $f_p: U_k\to Y$ with the property that $f_{p,1}: U_k\to Y$ is holomorphic and $K$-equivariant for all $p\in P$. Moreover, $f_{p,t}=f_p$ for $p\in P_0$, $t\in I$,  and $f_{p,t}=h$ on $X'\cap U_k$ for all $p\in P$, $t\in I$.  By restriction, we may assume that $U=U_k$ is an orbit convex Stein neighbourhood of $Kx_0$ such that $X_U=G\cdot U$ is Stein.

Let $R_U=R\cap U$ and let $\F$ denote the space of $K$-equivariant continuous maps from $(P\times R_U)\cup (P_0\times X_U)$ to $Y$ that are holomorphic along $X_U$ and equal to $h$ on $X'\cap X_U$.  Let $\F_U$ denote the corresponding space of maps with $X_U$ replaced by $U$.  By \cite[Lemma 2, p.~330]{HK1995}, the restriction map $\F\to\F_U$ is a homeomorphism, with each space given its compact-open topology.  Thus our homotopy $f_{p,t}$, which we constructed for elements of $\F_U$, lifts to a homotopy in $\F$.  Using Lemma \ref{l:homotopy-equiv}, this shows that $\Phi\hookrightarrow \Psi$ is a local weak homotopy equivalence.

It remains to construct our homotopy for $P$-families $\xi_p$ of $K$-equivariant continuous maps $U_k\to W$, which are holomorphic for $p\in P_0$ and vanish on $X'\cap U_k$. By Remark \ref{rem:P_0?}, we may assume that $\xi_p$ is holomorphic for $p$ in a neighbourhood $P_0'$ of $P_0$ in $P$. Let $\chi: P\to I$ be continuous such that $\chi=1$ on a neighbourhood of $P_0$ and $\chi=0$ on $P\setminus P_0'$. Let
\[  \tilde\xi_p=\chi(p)\xi_p,\quad \xi_{p,t}=t\tilde\xi_p+(1-t)\xi_p,\quad p\in P,\ t\in I.  \]
Clearly, $\tilde \xi_p$ is holomorphic for all $p\in P$, $\tilde \xi_p$ vanishes on $X'\cap U_k$, the homotopy is constant on a neighbourhood of $P_0$, and $\xi_{p,1}=\tilde\xi_p$ is holomorphic.
\end{proof}

Finally, we establish our approximation result.  It is an equivariant variant of \cite[Theorem 1]{Larusson2005}.  First, we prove a lemma.

\begin{lemma}  \label{lem:saturated}
Let $U$ be a $K$-invariant $\O(X)$-convex open subset of $X$.  Then $GU$ is saturated.
\end{lemma}

\begin{proof}
Since saturation is a local property over $Q$, we may assume that $\pi(U)$ is relatively compact.  Let $U_0$ be a relatively compact Stein neighbourhood of $\pi(U)$ in $Q$.  By Proposition \ref{prop:equivariant-embedding}(a), the Stein open subset $\pi^{-1}(U_0)$ of $X$ has a $G$-equivariant proper holomorphic embedding into a $G$-module.  Hence we may now assume that $X$ is a module.  For $x\in U$, let $F_x = \pi^{-1} (\pi(x))$ be the fibre through $x$.  By the Hilbert-Mumford criterion, there is a one-parameter subgroup $\lambda:\C^*\to G$, such that the closure of the orbit $\lambda(\C^*)x$ contains a point $y$ in the unique closed $G$-orbit $O$ in $F_x$. By \cite[Lemma 4.11]{Schwarz1989}, which applies since we are working in a module, we may assume that $\lambda(\mathbb T) \subset K$ (here $\mathbb T$ is the unit circle in $\C^*$). The holomorphic map $\eta:\C^*\to X$, $\theta \mapsto \lambda(\theta)x$, extends to $\C$ by setting $\eta(0) =y$.  Since $U$ is $K$-invariant and $\O(X)$-convex, the $\O(X)$-hull of $\lambda(\mathbb T)x$ is contained in $U$.  Since $0 \in \C$ is in the $\O(\C)$-hull of $\mathbb T$, we conclude that $y=\eta(0)\in U$, so $O\subset GU$.    All the $G$-orbits in $F_x$ have $O$ in their closure, so $GU$ intersects every $G$-orbit in $F_x$ and thus contains $F_x$. 
\end{proof}

\begin{proof}[Proof of Theorem \ref{t:main-theorem}(d)]  By assumption, the continuous $K$-map $f:X\to Y$ is holomorphic on a neighbourhood of a $K$-stable $\O(X)$-convex compact subset $L$ of $X$.  Take a relatively compact Stein neighbourhood $S_0$ of $\pi(L)$ in $Q$.  By Proposition \ref{prop:equivariant-embedding}(a), the saturated Stein open subset $S=\pi^{-1}(S_0)$ of $X$ has a $G$-equivariant proper holomorphic embedding in a $G$-module $V$.  For some of the arguments below, it helps to work in the smooth ambient space $V$ instead of the possibly singular space $X$.

In $V$, every neighbourhood of $L$ contains a smaller neighbourhood $U_0$ that is a sublevel set of a real-analytic $K$-invariant strictly plurisubharmonic exhaustion $\rho$ on $V$.  The existence of such exhaustions is well known, but for lack of a detailed reference, we will sketch a construction.  First, it is standard that $V$ carries a smooth strictly plurisubharmonic function $\rho_1$ which is negative on $L$ and positive outside any given neighbourhood of $L$ (see \cite[Proposition 2.5.1]{Forstneric2017}).  It is also standard that $\rho_1$ can be modified to an exhaustion $\rho_2$ of $V$ with those same properties.  Next, by density of real-analytic functions among smooth functions in the Whitney $C^2$ topology, $\rho_2$ may be approximated by a real-analytic function $\rho_3$ on $V$, closely enough that the properties of $\rho_2$ are preserved.  Finally, average $\rho_3$ over $K$ to obtain $\rho$.  For the reader's convenience, we insert here a proof that averaging preserves real-analyticity.  See \cite{Illman2003} for relevant results.

\begin{lemma}
Let the compact Lie group $K$ act continuously on a real-analytic manifold $M$.  The average over $K$ of a real-analytic function on $M$ is real-analytic. 
\end{lemma}

\begin{proof}
Note that the $K$-action on $M$ is automatically real-analytic.  By \cite[Sec.~6, Extension Theorem]{Heinzner1993}, there is a proper embedding of $M$ into a complex $G$-manifold $\widetilde M$.  Any real-analytic function $f$ on $M$ has a holomorphic extension $\tilde f$ to a neighbourhood $U$ of $M$ in $\widetilde M$.  We may take $U$ to be $K$-stable.  Holomorphic functions remain holomorphic when averaged.  Thus the average of $\tilde f\vert_U$ over $K$ is holomorphic and its restriction to $M$, which is the average of $f$ over $K$, is real-analytic.
\end{proof}

Choose a neighbourhood $U_0$ as above such that $f$ is holomorphic on $U=U_0\cap S$, viewed as a neighbourhood of $L$ in $X$, and let $\rho$ be a corresponding exhaustion.  Then $U$ is relatively compact, $K$-stable, and Stein, and moreover orbit convex by \cite[Lemma 1, p.~328, and Claim, p.~329]{HK1995}.  (For the definition and basic properties of orbit convexity, see \cite[Section 3.2]{Heinzner1991}.)  By \cite[Extension Lemma, p.~636, and Complexification Theorem (ii), p.~631]{Heinzner1991}, $U^\C=GU$ is a $G$-stable Stein neighbourhood of $L$ in $S$.  By \cite[Complexification Theorem (iii)]{Heinzner1991}, $f\vert_U$ extends to a holomorphic $G$-map $\tilde f:GU\to Y$.  By $G$-equivariance, we may extend $\tilde f$ further to a holomorphic $G$-map $GU\cup X'\to Y$ with $\tilde f=f$ on $X'$.
 
Recall that $U$ is a sublevel set of $\rho$ and let $R_\rho$ be the Kempf-Ness set in $V$ associated to $\rho$.  Then $R_\rho \cap U = R_\rho \cap GU$.  To see this, note that the intersection of $R_\rho$ with a fibre $F$ contains precisely the points where $\rho\vert_F$ attains its minimum; these points form a unique $K$-orbit contained in the unique closed $G$-orbit in $F$.  The $K$-equivariant deformation retraction of $V$ onto $R_\rho$ constructed in \cite{HH1994} preserves $GU$ since $GU$ is saturated by Lemma \ref{lem:saturated}.  It also preserves $U$ since $\rho$ decreases with the deformation.  Thus $U$ and $GU$ deformation-retract $K$-equivariantly onto the same subset $R_\rho \cap U = R_\rho \cap GU$ by a deformation $\alpha_t:GU\to GU$, $t\in [0,1]$, such that $\alpha_0$ is the identity and $\alpha_1$ is a retraction onto $R_\rho \cap GU$.  The deformation preserves $S\cap X'$ since $X'$ is closed and $G$-stable.

Now we shrink $U$ slightly and choose a continuous cut-off function $\theta$ on $Q$ with $\theta=1$ near $\pi(U)$ and $\theta=0$ outside a small neighbourhood of $\pi(U)$.  Define a continuous $K$-map $\beta:X\times[0,1]\to X$ by $\beta_t(x) = \alpha_{t\theta(\pi(x))}(x)$ for $x$ near $GU$ and $\beta_t(x)=x$ for $x$ away from $GU$.  Note that $\beta_t$ preserves $GU\cup X'$, $\beta_0=\id_X$, and $\beta_1(GU\cup X')\subset (R_\rho\cap GU)\cup X'$.  Therefore, on $GU\cup X'$, both $f$ and $\tilde f$ are $K$-homotopic to $f\circ\beta_1=\tilde f\circ\beta_1$ and thus $K$-homotopic to each other.  Extendability is deformation-invariant, so since $f\vert_{GU\cup X'}$ extends to a continuous $K$-map $X\to Y$ (namely $f$ itself), so does $\tilde f$.

Thus we may assume that $\tilde f$ extends to a continuous $K$-map $X\to Y$, also denoted $\tilde f$.  By Lemma \ref{lem:saturated} and Proposition \ref{prop:equivariant-embedding}(a), $GU=\pi^{-1}(\pi(U))$ has a $G$-equivariant proper holomorphic embedding $h$ in a $G$-module $W$.  By the equivariant Oka-Weil theorem (a special case of Theorem \ref{thm:C-O-W}), we can approximate $h:GU\to W$ uniformly on $L$ by a holomorphic $G$-map $k:X\to W$.  The inclusion $i:GU\hookrightarrow X$ factors $G$-equivariantly through the Stein $G$-space $M=X\times W$ (with the diagonal action) as $GU\stackrel{j}{\to} M \stackrel{\pi}{\to} X$, where $j=(i,h)$ is an embedding and $\pi$ is the projection.  The continuous $K$-map $\tilde f\circ\pi:M\to Y$ is holomorphic on $j(GU)\cup (X'\times W)$.  By Theorem \ref{t:main-theorem}(c), the restriction of $\tilde f\circ\pi$ to $j(GU)\cup (X'\times W)$ has a $G$-equivariant holomorphic extension $F:M\to Y$.  Then $F\circ (\id_X, k):X\to Y$ is a holomorphic $G$-map that approximates $f$ uniformly on $L$ and equals $f$ on $X'$.
\end{proof}

\begin{remark}   \label{rem:final-remark}
The following theorem would simplify the proof above by removing the need to move into a module.  For want of a detailed reference we have refrained from using it, but we indicate how a proof would proceed.

\smallskip

\noindent
{\bf Theorem.}  {\it Let $L$ be a $K$-stable $\O(X)$-convex compact subset in a Stein $G$-space $X$ and $U$ be a $K$-stable neighbourhood of $L$.  Then there exists a $K$-invariant real-analytic strictly plurisubharmonic exhaustion $\rho : X \to\R$ such that the corresponding Kempf-Ness set is a $K$-equivariant strong deformation retract of $X$ by a deformation that preserves the closure of each $G$-orbit in $X$ and such that
\[  L \subset \rho^{-1}(-\infty, 0) \Subset U \tag{$\ast$} \label{ast}.\]}

\smallskip

A proof of the theorem without condition \eqref{ast} was sketched in \cite{HK1995}.  There the authors explain how to build on ideas of Narasimhan \cite{Narasimhan1961} to make the deformation retraction work as in the smooth case, treated in \cite{HH1994}.  They ensure that the function $\rho$ is locally the restriction of a $K$-invariant strictly plurisubharmonic function defined on an open subset of a module, so the deformation is locally given as the restriction of a gradient flow from the smooth setting.  

A simple addition to this construction provides condition \eqref{ast}.  In the notation of \cite{HK1995}, the level sets $C_j$ of a strictly plurisubharmonic exhaustion can be replaced by any exhausting sequence of $K$-stable $\O(X)$-convex compact subsets $C_j$ with $C_j \subset \mathring{C}_{j+1}$.  For our purpose, we choose $C_0 = L$ and $C_1 \subset U$.  If, moreover, the $\epsilon_i$ are chosen so small that $\sum\epsilon_i < \frac 1{10} $, then the resulting real-analytic strictly plurisubharmonic exhaustion $\rho$ will have $\rho<\frac 1{10}$ on $L$ and $\rho\geq\frac 4 5$ on $X\setminus C_1$.  The function $\rho - \frac 1 2$ will then satisfy condition \eqref{ast}. 
\end{remark}

\end{document}